\documentclass{amsart}
\usepackage[utf8]{inputenc}

\title{Dehn--Sydler--Jessen Via Homological Algebra}
\author{Anubhav Nanavaty}

\usepackage{amsmath}
\usepackage{spectralsequences} 
\usepackage{graphicx}
\usepackage{scalerel}
\usepackage{stackengine,wasysym}
\usepackage{hyperref}

\usepackage{amssymb}
\usepackage{amsthm}
\usepackage{mathtools}
\usepackage{adjustbox}
\usepackage{dsfont}
\usepackage{MnSymbol}
\usepackage{mathrsfs}
\usepackage{stmaryrd}
\usepackage[all]{xy}
\usepackage[mathcal]{eucal}
\usepackage{verbatim}  
\usepackage{fullpage}  
\usepackage{hyperref}
\usepackage{setspace}
\usepackage{tikz-cd}
\usetikzlibrary{positioning}
\newcommand{\C}{\mathbb{C}}

\newcommand{\E}{\mathbb{E}}
\newcommand{\F}{\mathbb{F}}

\renewcommand{\H}{\mathbb{H}}

\newcommand{\N}{\mathbb{N}}

\newcommand{\Q}{\mathbb{Q}}
\newcommand{\R}{\mathbb{R}}

\newcommand{\Z}{\mathbb{Z}}

\renewcommand{\l}{\langle}
\renewcommand{\phi}{\varphi}
\renewcommand{\r}{\rangle}

\providecommand{\rm}[1]{\mathrm{#1}}

\usepackage{rotating}
\usepackage{amssymb}

\newtheorem{theorem}{Theorem}[section]

\newtheorem{proposition}[theorem]{Proposition}
\newtheorem{lemma}[theorem]{Lemma}
\newtheorem{corollary}[theorem]{Corollary}

\newtheorem{remark}[theorem]{Remark}

\theoremstyle{definition}
\newtheorem{definition}[theorem]{Definition}

\newtheorem{exercise}{Exercise}

\numberwithin{equation}{subsection}
\begin{document}
\maketitle
\begin{abstract}
    We provide an expository introduction to Euclidean Scissors Congruence, the study of polytopes in Euclidean space up to `cut and paste' relations. We first re-frame questions in scissors congruence as those in group homology. We then use this perspective to review the proof of the Dehn--Sydler--Jessen Theorem as found in the works of Dupont and Sah. 
\end{abstract}
\section{Introduction}
In Euclidean space, two polygons have the same area if and only if they are scissors
congruent: one can be cut up into finitely many smaller polygons and rearranged to equal the other. In
higher dimensions, understanding the relationship between volume and scissors congruence has proven more
complicated. Hilbert’s third problem asked: are two polyhedra scissors congruent if they have the same volume? His student, Max Dehn, constructed an invariant to answer this question in the negative \cite{Dehn},
which later culminated in the Dehn--Sydler--Jessen theorem \cite{Jessen} in the late 60’s. Two decades later, work
of Dupont and Sah \cite{DupontSahEuclidean} gave a different proof of the Dehn--Sydler--Jessen theorem via group homology computations. The purpose of this note is to give a self-contained review of this proof, parts of which appear in other work of Dupont and Sah (\cite{Dupont} and \cite{DupontFlags}). An exposition of the original proof is given in \cite{Schwartz}, and it is an interesting question to determine the relationship between these proofs. The statement of the theorem itself is as follows:
\begin{theorem}[Dehn--Sydler--Jessen]
There is an exact sequence of groups:
\begin{equation}\label{DehnSydlerJessen}
0\to \mathcal{P}(\E^2)\xrightarrow{-\times [0,1]} \mathcal{P}(\E^3)\xrightarrow{D} \R\otimes\R/\Z\xrightarrow{\phi}\Omega^1_{\R/\Z}\to 0
\end{equation}
where
\begin{enumerate}
    \item $\mathcal{P}(\E^n)$, the scissors congruence group of Euclidean $n$ space $\E^n$
    \item The map $-\times [0,1]$ is induced by the inclusion of a polygon $P\subset\E^2$ as a prism $P\times [0,1]\subset \E^3$
    \item The map $D$ is the \emph{Dehn invariant}
    \item The group $\Omega^1_{\R}$ of K\"{a}hler differentials of $\R$ over $\Z$.
    \item The map $\phi$ sends $\ell\otimes \theta/\pi$ to $\ell\frac{d\cos(\theta)}{\sin(\theta)}$
    \end{enumerate}
\end{theorem}
The second section of these notes defines \emph{the polytope module}, $\mathrm{Pt}(\E^n)$, and explains how to rephrase the computations of $\mathcal{P}(\E^n)$ to those of group homology with coefficients in $\mathrm{Pt}(\E^n)$. The key theorem that makes the homological algebra approach work is the following:
\begin{theorem}\label{thm:TitsComplex}
As $\mathrm{Isom}(\E^n)$-modules, we have an isomorphism:
    \[\tilde{H}_{*}(\mathscr{T}(\E^n),\Z)=\begin{cases}Pt(\E^n)\quad&*=n-1 \\
    0\quad&\text{ otherwise}
    \end{cases}\] where $\mathscr{T}(\E^n)$ is the \emph{Tits complex} of $\E^n$, a simplicial set defined by:
    \[\mathscr{T}(\E^n)_k:=\Big\{( U_0\supset\dots \supset U_k):U_i\text{ are affine subspaces of $\E^n$ and } U_0,U_k\neq \O,\E^n\Big\}\]
    with face and degeneracy maps given by the usual deletion and insertion maps respectively.
\end{theorem}
While this review closely follows \cite{Dupont}, one should also read this section alongside Section 2 of \cite{Malkiewich} which provides spectral lifts of these constructions.

In the third section, we study what happens if we classify polytopes in $\E^n$ only up to translation, reducing the problem to understanding the $O(n)$ orbits. The last two sections focus on the 3-dimensional setting. Section 4 deals constructs the exact sequence (see \ref{thm:almostDSJ}):
\[
    0\to H_2(O(3),(\E^3)^t)\to H_0(O(3),\mathcal{D}^1(\E^3))\to \R\otimes\R/\Z\to H_1(O(3),(\E^3)^t)\to 0
    \]
where the map $H_0(O(3),\mathcal{D}^1(\E^3))\to \R\otimes\R/\Z$ recovers the Dehn invariant. In the last section, we put together theorems from the previous sections to prove Dehn--Sydler--Jessen, along with computations in group homology special to low degree:
\begin{theorem}\label{thm:exceptional}
    \[H_1(O(3),(\E^3)^t)=\Omega_{\R}^1\]
    and 
    \[H_2(O(3),(\E^3)^t)=0\]
\end{theorem}
These proofs rely heavily on the structure of the quaternion algebra $\H$. Considering $\E^3$ as a subspace of $\H$, all of a sudden previously intractable computations become possible. The appendix finishes a particular computation in Section 5 that detracts from the narrative. 
\\\\\\
\textbf{Acknowledgements.} The author would like to thank Cary Malkiewich, Daniil Rudenko, and Jesse Wolfson for helpful discussions. 
\tableofcontents
\section{The Polytope Module and the Tits Complex}
We begin with some definitions:
\begin{definition}
   Let $\E^n$ denote Euclidean space, i.e. $\R^n$ equipped with the standard inner product. Let $\mathrm{Isom}(\E^n)$ denote the isometry group of $\E^n$.
\end{definition}
The main objects of study here are polytopes, and we must make precise what we mean by them:
\begin{definition}
    A \emph{geodesic} $n$-simplex in $\R^n$ is the convex hull of $n+1$ points $\{a_0,\dots,a_n\}$, which we write $|(a_0,\dots,a_n)|$. A \emph{Polytope} $P$ is a finite union:
    \[
    P=\bigcup_{i=1}^n\Delta_i
    \]
    where $\Delta_i$ are geodesic $n$-simplies, and the intersection $\Delta_i\cap \Delta_j$ is a geodesic $k$-simplex for $k<n$. Lastly, if $P,P',P''$ are polytopes, then we write 
    \[P=P'\coprod_\circ P''\] 
    If $P=P'\cup P''$ and $P'\cap P''$ has no interior points.
\end{definition}

\begin{definition}
     Let $\mathrm{Isom}(\E^n)$ denote the isometry group of $\E^n$. The \emph{Polytope module} is the left $\mathrm{Isom}(\E^n)$-module defined as
    \[\mathrm{Pt}(\E^n):=\Z\Big\l [P]: [P]=[P'']+[P']\text{ for }P=P'\coprod_\circ P''\Big\r\]
    where $[P]$ ranges over all polytopes. 
\end{definition}
Classical scissors congruence comprises of the study of the module structure of $\mathrm{Pt}(\E^n)$. The purpose of this section is to show that one can understand $\mathrm{Pt}(\E^n)$ using homological methods.
\begin{definition}\label{simplicial chains}
    Let $C_*(\E^n)$ denote the chain complex with 
    $C_n(\E^n)=\Z\l (a_0,\dots,a_n):a_n\in\E^n \r$ and the boundary maps generated by
    \[
    \partial_n(a_0,\dots,a_n)=\sum_{i=1}^n(-1)^i(a_0,\dots,\hat{a_i},\dots,a_n)
    \]
    There is an increasing filtration of chain complexes
    \[
    \O=\Big(C_*(\E^n)\Big)_{-1}\subset \Big(C_*(\E^n)\Big)_{0}\subset\dots\subset \Big(C_*(\E^n)\Big)_{n}=C_*(\E^n) 
    \]
    where $(a_0,\dots,a_k)\in C_k(\E^n)_p$ if its convex hull $|(a_0,\dots,a_k)|$ is contained in a geodesic $p$-simplex. We can similarly define $C_*(V)$ for any geodesic subspace $V\subset \E^n$.
\end{definition}
\begin{definition}\label{twist}
    For the purposes of this note, we view $C_n(\E^n)$ as an $Isom(\E^n)$-module with action given by $\phi\cdot (a_0,\dots,a_n)=\mathrm{sgn}(\phi)(\phi(a_0),\dots,\phi(a_n))$, where $\mathrm{sgn}(\phi)=\pm 1$ is $1$ iff $\mathrm{det}(\phi)>0$ and $-1$ if $\mathrm{det}(\phi)<0$. To recall this fact, we write $C_n(\E^n)^t$, where $t$ stands for `twist'
\end{definition}
The reason for this will become apparent in the next proposition:
\begin{proposition}\label{prop:chains}
We have a map of $\mathrm{Isom}(\E^n)$-modules:
\[
\phi:C_n(\E^n)\to \mathrm{Pt}(\E^n)
\]
given by fixing an orientation on $\E^n$ and defining, for $\sigma=(a_0,\dots,a_n)$: \[\phi(\sigma)=\begin{cases} \epsilon_\sigma |\sigma|\quad\quad&\sigma\in \mathrm{Gr}_n(C_n(\E^n))\\
0\quad\quad&\text{otherwise}
\end{cases}\]
where $\mathrm{Gr}_n(C_n(\E^n))$ denotes the $n$-th graded piece of the filtration on $C_n(\E^n)$, and $\epsilon_\sigma=\pm 1$ and is positive iff the orientation of $\sigma$ is the same as the standard orientation on $\E^n$ given by the ordering of the vertices.  
This induces an isomorphism:
    \[
    H_n\Big(Gr_n(C_*(\E^n))\Big)\xrightarrow{\cong }\mathrm{Pt}(\E^n)
    \]

\end{proposition}
\begin{proof}
    We first prove that the map $\phi$ on $Gr_n(C_n(\E^n))$ descends to a map on homology of this graded complex. Since the boundary map induced on the graded pieces is just the zero map, we have:
    \[
    H_n(Gr_n(C_*(\E^n)))=C_*(\E^n)/[C_n(\E^n)_{n-1}+\partial C_{n+1}(\E^n)]
    \]
We already know that if $\sigma\in C_n(\E^n)_{n-1}$, then $\phi(\sigma)=0$. Taking $\sigma\in C_{n+1}(\E^n)$, we show $\phi(\partial \sigma)=0$ as follows:
\begin{enumerate}
    \item $\sigma=(a_0,\dots,a_{n+1})$ defines a convex hull of $n+1$ points inside $\E^n$, so for each $i$, its faces $\sigma_{\hat{i}}:=(a_0,\dots,\hat{a_i},\dots, a_{n+1})$ are either non-zero in $Gr_n(C_{n}(\E^n))$ or are in $(C_{n}(\E^n))_{n-1}$
    \item $\sigma_{\hat{i}}\neq 0$ in $Gr_n(C_{n}(\E^n))$, then we have     $\phi(\sigma_{\hat{i}})=\epsilon_{i}|\sigma_{\hat{i}}|$, where $\epsilon=\pm 1$ depending on whether the ordering orientation on $\sigma$ agrees with  that of $\E^n$, 
    and so \[\phi(\partial\sigma)=\sum_{\sigma_{\hat{i}}\text{ proper}}(-1)^i\epsilon_{i}|\sigma_{\hat{i}}|\]
    \item For each tuple where $0\leq i<j\leq n$, there is a unique supporting hyperplane going through points $(a_1,\dots,\hat{a_i},\dots,\hat{a_j},\dots,a_n)$ splitting $\sigma$ either into two simplices or none at all, which we write as $P_{i,j}$. Note that $P_{i,j}$ will always intersect $\sigma_{\hat{i},\hat{j}}:=|(a_0,\dots,\hat{a_i},\dots,\hat{a_j},\dots,a_n)|$, but may give a degenerate simplex. Write $P_{i,j}^1,P_{i,j}^2$ as the different sub-simplices of $\sigma$ that $P_{i,j}$ divides (if it does). Then we can write the sum in $\mathrm{Pt}(\E^n)$ as: \[|\sigma_{\hat{i}}|=\sum_{j\neq i}|\sigma_{\hat{i}}\cap P_{i,j}^{1}|+|\sigma_{\hat{i}}\cap P_{i,j}^{2}|\]
    as the other possible hyperplanes will \emph{definitely} not cut $\sigma_{\hat{i}}$. 
    \item By 3, it follows that that $|\sigma_{\hat{i}}\cap P_{i,j}^{1}|$ will appear in the sum $\sigma_{\hat{j}}$ as well! We now show that it appears with opposite sign. First, the sign on $|\sigma_{\hat{i}}\cap P_{i,j}^{1}|$ is $(-1)^i\epsilon_i$, where as the sign for $|\sigma_{\hat{j}}\cap P_{i,j}^{1}|$ is $(-1)^j\epsilon_j$. If $i,j$ have the same parity, then the number of transpositions taking $(1,\dots,\hat{i},\dots,n)$ to $(1,\dots,\hat{j},\dots,n)$ is odd, meaning they have different orientation, i.e. $\epsilon_i=-\epsilon_j$. A similar argument shows that if $i$ and $j$ have different parity, then $\epsilon_i=\epsilon_j$. The same argument works for $|\sigma_{\hat{i}}\cap P_{i,j}^{2}|$.  
    \item Since each $|\sigma_{\hat{i}}\cap P_{i,j}^{1}|$ and $|\sigma_{\hat{i}}\cap P_{i,j}^{2}|$ in $\mathrm{Pt}(\E^n)$ appears not at all, or twice with different parity in the sum of $\phi(\partial\sigma)$, we conclude that $\phi(\partial\sigma)=0$. 
    \end{enumerate}
So, we have that the map $\phi$ is well-defined on homology, and one can see that it is surjective just by writing a geodesic simplex $\Delta^n$ in $\E^n$ as a convex hull of points $\sigma=(a_0,\dots,a_n)$, and so $\phi(\sigma)=\Delta^n$ (up to sign). Showing injectivity is a little more painful. We construct an inverse map
\[
\psi:\mathrm{Pt}(\E^n)\to H_n(Gr_n(C_*(\E^n)))
\]
as follows: we know that $\mathrm{Pt}(\E^n)$ is generated by geodesic $n$-simplicies. For such a geodesic $n$-simplex $\Delta^n$, we can write it as the convex hull of points $|(a_0,\dots,a_n)|$ labeled in such a way that agrees with the orientation on $\E^n$, then define $\psi(\Delta^n)=(a_0,\dots,a_n)$. In order to ensure this map is well-defined, we need to make sure that $\psi$ agrees if we triangulate $\Delta^n$ into a different simplicial complex $|K|$, then $\psi(|K|)=\psi(\Delta^n)$. By the simplicial approximation theorem, we can barycentrically subdivide $Sd^{n}(K)$ and get a simplicial map $f:Sd^{n}(K)\to \Delta^n$ that is an approximation of the identity map $id:\Delta^n\to \Delta^n$, and since $f$ is simplicial, $f$ maps $\partial (Sd^{n}(K))$ to $\partial \Delta^n$. This means that $f_*:C_*(\Delta^n)\to C_*(\Delta^n)$ is chain homotopic to the identity, i.e. by definition there are morphisms $h_n:C_n(\Delta^n)\to C_{n+1}(\Delta^n)$ such that
\[f_*[Sd^{n}(K)]-[Sd^{n}(K)]=\partial h_n[Sd^{n}(K)]+h_{n-1}\partial[Sd^{n}(K)]\]
Since we are working with a nice simplicial map, we can construct $h_n$:
\begin{equation}\label{$h_n$ formula}
h_n:\sum_{\sigma\in\{\text{n-simplices of $K$}\}}\sum_{i=0}^n(-1)^i\big(a_0^\sigma,\dots,a_i^\sigma,f(a_{i}^\sigma),\dots,f(a_n^\sigma)\big)
\end{equation}
We also see that if $f_*$ sends $\partial[Sd^nK]$ to $\partial\Delta^n$, we also have that on relative homology $f_*[Sd^{n}(K)]=[\Delta^n]$, which implies that in $C_n(\Delta^n)$, 
\begin{equation}\label{rel hom formula}
f_*[Sd^{n}(K)]=[\Delta^n]+Q
\end{equation} where $Q$ is a sum of degenerate simplicies that live in $\partial\Delta^n$. Putting this together, we have in $C_n(\Delta^n)$:
\[
[\Delta^n]-[Sd^{n}(K)]=\partial h_n[Sd^{n}(K)]+h_{n-1}\partial[Sd^{n}(K)]-Q
\]
By looking at the formula for $h_n$ in \ref{$h_n$ formula}, we see that $h_{n-1}\partial[Sd^{n}(K)]\in \Big(C_*(\Delta^n)\Big)_{n-1}$ which will also hold if we look at $C_*(\E^n)$, and combing this with \ref{rel hom formula}:
\[
[\Delta^n]-[Sd^{n}(K)]\in \partial_{n+1}(C_{n+1}(\E^n))+(C_{n}(\E^n))_{n-1}
\]
A similar, but easier argument shows that:
\[
\psi(K)-\psi(Sd^{n}(K))\in (C_{n}(\E^n))_{n-1}
\]
And so we have $[\Delta^n]-[K]=0$ in $H_n(\mathrm{Gr}_n(C_*(\E^n))$, meaning that $\psi[\Delta^n]=\psi[K]$, so the map is well-defined. One can see that $(\psi\circ \phi)[K](a_0,\dots,a_n)=(a_0,\dots,a_n)$, and so $\phi$ is a surjective homomorphism with an inverse, making it an isomorphism. Lastly, one can check that both $\phi$ and $\psi$ respect the action of $\mathrm{Isom}(\E^n)$, which allows us to promote the isomorphism to one of $\mathrm{Isom}(\E^n)$-modules.
\end{proof}
The main purpose of the last proposition is that it is easier to compare $C_*(\E^n)$ with homology of the \emph{Tits} complex $\mathcal{T}(\E^n)$: 
\begin{theorem}[Theorem \ref{thm:TitsComplex}]
   As $\mathrm{Isom}(\E^n)$-modules, we have an isomorphism:
    \[\tilde{H}_*(\mathcal{T}(\E^n),\Z)\cong \begin{cases}
        \mathrm{Pt}(\E^n)\quad&*=n-1\\
        0\quad&\text{otherwise}
    \end{cases}\]
\end{theorem}
\begin{proof}
    Consider the following double complex, where:
    \[
A_{p,*}=\begin{cases}C_*(\E^n)_{n-1}\quad &\text{$p=-1$}\\
\bigoplus_{\Phi=(U_0\supseteq\dots\supseteq U_p)}C_*(U_p)\quad&\text{$p>-1$}\\
\end{cases}
\]
Where the horizontal maps $\partial'_{p,q}:A_{p,q}\to A_{p+1,q}$ are given by the face maps of the Tits complex, i.e. $\sum_{i=0}^p(-1)^is^i$, and the vertical maps are given as $\partial''_{p,q}=(-1)^p\partial_q$, where $\partial_q:C_q(U_p)\to C_{q-1}(U_p)$ is the boundary map on simplicial chains as in \ref{simplicial chains}. When $p=0$, the maps $\bigoplus_{U_0}C_*(U_0)\to C_*(\E^n)_{n-1}$ are simply the inclusion maps. We show that the  horizontal boundary maps are exact - consider the map:
\begin{equation}\label{nullhomotopy}
s_p:A_{*,p}\to A_{*,p+1};\quad s_p(\sigma_{(U_0\supset\dots\supset U_p)})=\sigma_{(U_0\supset\dots\supset U_p\supset U_\sigma)}
\end{equation}
where $U_\sigma=\mathrm{Span}(\sigma)$. One can check that this defines a null-homotopy, and so all of the homology of the total complex will lie in degree $-1$: 
\[
H_*\Big(\mathrm{Tot}(A_{*,*})\Big)=H_*(A_{-1,*},\partial')=H_*(C_*(\E^n)_{n-1},\partial')
\]
Given the short exact sequence of chain complexes:
\[
0\to C_*(\E^n)_{n-1}\to C_*(\E^n)_{n}\to Gr_n(C_*(\E^n))\to 0
\]
the long exact sequence in homology tells us that
\[
\tilde{H}_*(C_*(\E^n)_{n-1})=H_{*+1}(Gr_n(C_*(\E^n)))
\]
Now, looking at the vertical maps in $A_{p,*}$, we see that they compute the homology of the complex $C_*(U_p)$, which is trivial except at $q=0$, where it is $\Z$, and so we also have that 
\[
H_*(\mathrm{Tot}(A_{*,*}))=H_*(A_{*,0},\partial''_{*,0})
\]
further, $A_{*,0}=\bigoplus_{(U_0\supset\dots \supset U_p)}\Z$, i.e. these chains compute the cohomology of $\mathcal{T}(\E^n)$. So, we have that \[H_*(\mathrm{Tot}(A_{*,*}))=H_*(\mathcal{T}(\E^n))\]
So, we can invoke \ref{prop:chains} conclude the theorem.
\end{proof}

    \begin{corollary}\label{lustigex}

    There is an exact sequence of $\rm{Isom}(\E^n)$-modules:
        \[
    0\to \mathrm{Pt}(\E^n)\xrightarrow{\partial_1} \bigoplus_{U_{n-1}}\mathrm{Pt}(U_{n-1})\xrightarrow{\partial_2}\cdots\xrightarrow{\partial_{n-2}} \bigoplus_{U_{2}}\mathrm{Pt}(U_{2})\xrightarrow{\partial_{n-1}}\bigoplus_{U_{0}}\mathrm{Pt}(U_0)\xrightarrow{\epsilon} \Z\to 0
    \]
    called the \emph{Lusztig exact sequence}.
    \end{corollary}
\begin{proof}
    We have a filtration of chain complexes on $C_*(V)$ for any affine subspace $V\subset\E^n$:
    \[
    \O=\Big(C_*(V)\Big)_{-1}\subset\Big(C_*(V)\Big)_{0}\subset \dots\subset \Big(C_*(V)\Big)_{n}=C_*(V)
    \]
    where $\Big(C_*(V)\Big)_{p}$ is defined as in \ref{simplicial chains}. Observe that if we look at the graded pieces of the filtration we have:
    \[
    Gr_k((C_*(V)))=\begin{cases}
        \bigoplus_{\mathrm{dim}(U_k)=k}C_k(U_k)/[C_k(U_k)]_{k-1}\quad& *=k\\
        0\quad&\text{otherwise}
    \end{cases}
    \]
    where $U_k$ ranges over all the affine subspaces. This filtration gives us a spectral sequence with $E^1$ page
   \[
   E^1_{p,q}=H_{p+q}\Big(Gr_p((C_*(V)))\Big)
   \]
   Using our knowledge of $Gr_k((C_*(V)))$ and following the proof of \ref{prop:chains}, we have:
   \[
   E^1_{p,q}=\begin{cases}
       \bigoplus_{U_p}\mathrm{Pt}(U_p)\quad&q=0\\
       0\quad&q\neq0
   \end{cases}
   \]
   so, $E^2=E^\infty$. Now, the differential $E^1_{p,q}\to E^1_{p-1,q}$ is induced by the boundary map $C_q(U_p)\to C_{q-1}(U_p)$, which we know is exact except with homology $\Z$ in degree $0$, meaning that $E^2_{p,q}=0$ for $(p,q)\neq (0,0)$ and $E^{2}_{0,0}=\Z$. So, we have an exact sequence:
\[
0\to E^1_{n,0}\to E^1_{n-1,0}\to\dots\to E^1_{0,0}\to \Z\to 0
\]
which is what we wanted. 
\end{proof}
\section{Translational Scissors Congruence}
Now, we can define $G$-scissors congruence, and rephrase computations as those coming from group homology.
\begin{definition}
    Let $G<\mathrm{Isom}(\E^n)$. Then the $G$-scissors congruence group is defined as 
    \[
    \mathcal{P}(\E^n,G):= \mathrm{Pt}(\E^n)/\l [P]-[\phi(P)]:\phi\in G\r
    \]
    where $\phi(P)$ is the image of $P$ under the map $\phi:\E^n\to \E^n$. Define \emph{the} scissors congruence group of $\E^n$ to be:
    \[
    \mathcal{P}(\E^n):=\mathcal{P}(\E^n,\mathrm{Isom}(\E^n))
    \]
\end{definition}
The following we leave as an exercise:
\begin{exercise}
    Let $T(n)$ denote the group of translations in $\E^n$, and $O(n)$ denote the orthogonal group. Then:
    \[
    \mathrm{Isom}(\E^n)=T(n)\rtimes O(n)
    \]
\end{exercise}
One can phrase these questions in terms of \emph{group} homology, defined below:
\begin{definition}
    Let $\Z[G]$ denote the group ring of $G$, the free abelian group with generators corresponding to $g\in G$, and multiplication corresponding to the group operation in $G$. By $G$-module, we mean a module over the ring $\Z[G]$.
\end{definition}
\begin{definition}
    Let $C_*(G)$ be the same chain complex as in \ref{simplicial chains}, but where $C_n(G)=\Z[G^n]$. For any $G$ module $M$, we define $C_*(G;M):=C_*(G)\otimes M$. For each $n\in\N$, $C_n(G;M)$ is a left $G$ module with action given by
    \[
    g\cdot (x_1,\dots,x_n,m)=(gx_1,\dots,gx_n,gm)
    \]
    and the usual differential $\partial:C_n(G)\to C_{n-1}(G)$ extends to a differential on $C_n(M)$, and commutes with the action of $G$. Let $(C_*(G;M))_G$ denote the complex of coinvariants of $C_*(G;M)$, i.e. where we quotient $C_*(G;M)$ by the $G$-action. We define the \emph{group homology of $G$ with coefficients in $M$ as }:
    \[
    H_*(G;M)=H_*(C_*(G;M))_G)
    \]
    and when $M=\Z$ with trivial action of $G$, we write $H_*(G)=H_*(C_*(G;\Z)_{G})$.
\end{definition}
We can now make the following observation:
\begin{proposition}\label{prop:gpcoh}
    Let $N\unlhd \mathrm{Isom}(\E^n)$. Then, we have:
    \[
    \mathcal{P}(\E^n,G):=H_0\Big(G/N,H_0(N, \mathrm{Pt}(\E^n))\Big)
    \]
\end{proposition}
\begin{proof}
    First recall that for a left $G$-module $M$, we have
    \[
    H_0(G,M)=M_G
    \]
    where $M_G:=M/\{gm-m:g\in G,m\in M\}$ is the coinvariants of $M$. So, 
    \begin{align*}
    H_0(G/N,H_0(N, \mathrm{Pt}(\E^n)))&=H_0(G/N,\mathrm{Pt}(\E^n)_{N})\\
    &=(\mathrm{Pt}(\E^n)_{N})_{G/N}\\
    &=\mathrm{Pt}(\E^n)_{G}\\
    &=\mathcal{P}(\E^n,G)
    \end{align*}
\end{proof}

Proposition \ref{prop:gpcoh} immediately gives us the following corollary:
\begin{corollary}
\[
    \mathcal{P}(\E^n)=H_0\Big(O(n),H_0(T(n), \mathrm{Pt}(\E^n))\Big)
    \]
\end{corollary}
From this corollary, a good first question to ask is: how to compute group homology of $T(n)$? Here, we can use that $T(n)$ is a torsion-free abelian group:
\begin{theorem}\label{thm:abtorfree}
    If $V$ is a torsion-free abelian group, then addition $+:V\times V\to V$ endows $H_*(V)$ with the structure of an anti-commutative graded ring. Further, we have an isomorphism of graded rings:
    \[
    H_*(V)\cong \bigwedge\nolimits^*_\Z (V)
    \]
\end{theorem}
\begin{proof}
    For any group $G$, $H_0(G)=\Z[G_{G}]=\Z$, and so we move on to the general statement for $n>1$. We see that $C_1(V)_V=\Z[V^2]/\l(v_0,v_1)-\tau (v_0,v_1), \tau\in V\r$, so translating $v_0$ to the origin, we have the isomorphism $C_1(V)\cong \Z[V]$. Further, the kernel of $\partial_1$ will be exactly $V\subset \Z[V]$. We have for any $(v_0,v_1,v_2)\in C_2(V)$:
    \[\partial_2(v_0,v_1,v_2)=(v_0,v_1)-(v_0,v_2)+(v_1,v_2)=(v_1,v_1)\]
    which is the same as $(0,0)$ when taking coinvariants $C_1(V)_{V}$, so the functorial map $V\to C_1(V)$ sending $v\mapsto (0,v)$ gives an isomorphism
    \[
    V\xrightarrow{\cong} ker(\partial_1)/im(\partial_2)=H_1(V)
    \]
    By the K\"{u}nneth theorem we have for each $k\in\N$ an exact sequence of abelian groups:
    \[
    0\to \bigoplus_{i+j=k}H_i(V)\otimes H_j(V)\xrightarrow{\times} H_k(V^2)\to \bigoplus_{i+j=k-1}\mathrm{Tor}_1^\Z(H_i(V),H_j(V))\to 0
    \]
    In our case there is no torsion by assumption and so $\times$ is an isomorphism. Since $V$ is an abelian group, the addition map $+:V^2\to V$ is in fact a group homomorphism, and by functoriality we get a map $H_k(V^2)\to H_k(V)$ for each $k\in\N$. The \emph{Pontryagin product} is given as the composition of $+$ with $\times$: 
    \[   \bigwedge\nolimits_{i,j}:\bigoplus_{i+j=k}H_i(V)\otimes H_j(V)\to H_k(V)
    \]
    One can check that $\bigwedge\nolimits_{i,j}\Big((v_1,\dots,v_i),(v_1,\dots,v_j)\Big)=-\bigwedge\nolimits_{j,i}\Big((v_1,\dots,v_j),(v_1,\dots,v_i)\Big)$, making $H_*(V)$ a anti-commutative graded ring. By universal property of the exterior algebra over $\Z$, the functorial isomorphism $V\to H_1(V)=V$ extends to a functorial map of graded algebras:
    \[
    \psi:\bigwedge\nolimits^*_\Z (V)\to H_*(V)
    \]
    We prove that $\psi$ is an isomorphism for $V=\Z^n$, and the general case follows as any torsion-free abelian group can be presented as a direct limit of finitely generated free-abelian groups, and $H_*$ commutes with direct limits.
    We already have the identity $\bigwedge^k(V^2)\cong \bigoplus_{i+j=k}\bigwedge^i(V)\otimes\bigwedge^j(V)$, giving us a commutative diagram
    \[
    \begin{tikzcd}
         \bigoplus_{i+j=k}H_i(V)\otimes H_j(V)\arrow[r,"\cong"]& H_k(V^2)\\
        \bigoplus_{i+j=k}\bigwedge^i(V)\otimes\bigwedge^j(V)\arrow[r,"\cong"]\arrow[u]&\bigwedge^k(V^2)\arrow[u]
    \end{tikzcd}
    \]
    so we just need to prove that the vertical arrows are isomorphisms when $V=\Z$. We already have that $H_1(\Z)=\Z$, and $H_0(\Z)=\Z$, and so we must show $H_i(\Z)=0$ for $i>1$. This follows from the observation that one can write for $i>1$ an element of $C_i(\Z)_{\Z}$ uniquely as $(0,v_1,\dots,v_i)$, and so homology will be a acyclic by constructing a null-homotopy exactly as in \ref{nullhomotopy}. We have $\bigwedge\nolimits_\Z^*\Z=\Z\oplus\Z$ with a copy concentrated in degree $0$ and $1$ each, and by previous work, the induced map to $H_*(\Z)$ is an isomorphism, so we get the full result.
\end{proof}
The goal now is to use this fact to compute $\mathcal{P}(\E^n,T(n))$, which comes from the framework built in the previous section:
\begin{theorem}\label{thm:sseqtr}
Similar to the proof of \ref{thm:TitsComplex}, define the double complex
\[\tilde{A}_{p,*}=\begin{cases}
C_*(\E^n)_{T(n)}\quad& q=-1\\
    \bigoplus_{U_p\supset\dots\supset U_0}\ C_*(U_p)_{T(U_p)}\quad& q>-1\\
\end{cases}
\]
where $T(U_p)$ is the translation group of $U_p$, and $U_p$ is now a vector subspace of $\E^n$. The spectral sequence associated to the double complex $A_{p,q}$ has an $E^2$ page with
    \[
    E^2_{p,q}=\begin{cases}
        \tilde{H}_{p}(\mathscr{T}(V),\bigwedge\nolimits^q_\R(\mathfrak{g}))\quad\quad& p+q=n-1\\
        0\quad\quad& \text{otherwise}
    \end{cases}
    \]
    where $\bigwedge\nolimits^q_\R(\mathfrak{g})$ is the sheaf on 
    $\mathcal{T}(\E^n)$ sending $\bigwedge\nolimits^q_\R(\mathfrak{g}):(U_0\supset\dots\supset U_p)\mapsto \bigwedge\nolimits^q_\R(U_p)$ and $\bigwedge\nolimits^0_\R(\mathfrak{g})$ is the constant sheaf $\Z$. Further, $E^r$ differentials vanish for $r\geq 2$, and $H_{n-1}(\mathrm{Tot} \tilde{A}_{\bullet,\bullet})\cong 
   \mathcal{P}(\E^n,T(n))$.
\end{theorem}
\begin{proof}
Observe that $\partial:\tilde{A}_{p,q}\to \tilde{A}_{p,q-1}$ is the $E_0$ differential given by $(-1)^p$ times the usual differential of the bar complex. Using \ref{thm:abtorfree}, we see that taking homology with respect to this differential gives us:
\[
E_1^{p,q}=\begin{cases}
\Z\quad\quad& p=-1,\ q=0\\
\bigwedge\nolimits^q_\Z(\E^n)\quad\quad &p=-1,\ q\neq 0\\
\bigoplus_{\Phi=(U_0\supseteq\dots\supseteq U_p)}\bigwedge\nolimits^q_\Z(U_p)\quad\quad &p>-1
\end{cases}
\]
Since all the $U_i$ are $\R$-vector spaces, have a group isomorphisms for $n\geq 1$: $\bigwedge\nolimits^n_\Z (U_i)\cong \bigwedge\nolimits^n_\R (U_i)$. For dimension reasons, we see that $E^1_{p,q}=0$ for $p+q>n-1$, and so the same must hold for the $E^2$ page

As for the vanishing of the differentials, consider the \emph{dilation operators} $\{\mu_\lambda:\R^n\to\R^n:\lambda\in\R^\times\}$ such that $\mu_\lambda(x)=\lambda x$. Then we observe that since $\mu_\lambda$ commutes with the chain maps of $\tilde{A}_{p,q}$, it follows that \[\mu_\lambda^{p-2,q-1}\circ d^2_{p,q}=d^2_{p,q}\circ \mu_\lambda^{p,q}\]
    however, in $E^2_{p,q}$ we have that $\mu_\lambda^{p,q}$ is scaling by $\lambda^q$, and on $E^2_{p-2,q-1}$, it follows that $\mu_\lambda^{p-2,q-1}$ scales by $\lambda^{q-1}$. This implies that $\lambda^{q-1}d^2_{p,q}=d^2_{p,q}\lambda^{q}$, i.e. $\lambda d^2_{p,q}=1$ for all $\lambda\in\R^\times$, implying $d^2_{p,q}=0$. A similar argument shows that in fact $d^r=0$ for all $r\geq 2$. Lastly, for the same reason as in the proof of \ref{thm:TitsComplex}, the horizontal boundary maps are exact, giving us \[H_{n-1}(\mathrm{Tot}(\tilde{A}_{\bullet,\bullet}))=H_n(Gr_n(C_*(\E^n)_{T(n)}))=\mathrm{Pt}(\E^n)_{T(n)}\] which is equivalent to $\mathcal{P}(\E^n,T(n))$ by \ref{prop:gpcoh}. Since $d^r=0$ for $r\geq 2$, we have that for all $k\in \{0,\dots, n-2\}$
    \[
    H_n(Gr_n(C_*(\E^n)_{T(n)}))=\bigoplus_{p+q=n-1}H_{k+1}(Gr_n(C_*(\E^n)_{T(n)}))=0
    \]
    the latter being equal to zero because the $k$-chains of $Gr_n(C_*(\E^n)_{T(n)})$ are simply zero (every $k$-simplex is contained in an $n-1$-simplex). Therefore, $E^2_{p+q}=0$ for $p+q<n-1$ as well.
\end{proof}
We now have a greater understanding of the full scissors congruence group by studying the $E^\infty$ page of \ref{thm:sseqtr}. For $n=0$, we get $\mathcal{P}(\E^0)=\Z$ generated by the point, and for $n>0$ we have the decomposition:
\begin{corollary}\label{cor:decomp}
Write $\mathcal{D}^q=\tilde{H}_{n-q-1}(\mathscr{T}(\E^n),\bigwedge\nolimits^q_\R(\mathfrak{g}))$ where $0\leq q<n$ and write $\mathcal{D}^n:=\bigwedge\nolimits^n_{\R}\E^n$.  We then have:
    \[\mathcal{P}(\E^n,T(n))\cong \bigoplus_{q=1}^n H_0(O(n),\mathcal{D}^q)\]
\end{corollary}
\begin{remark}
    The geometric interpretation to $\mathcal{D}^q(\E^n)$ was first given by Hadwiger in \cite{Hadwiger}. Suppose $q=n$, then we have that if $\sigma=(a_0,\dots, a_n)$ is a non-degenerate simplex, the projection $\mathcal{P}(\E^n,T(n))\to \bigwedge\nolimits^n_{\R}\E^n$ will send $\sigma$ to $a_0\wedge\dots\wedge a_n$. Picking the isomorphism $\bigwedge\nolimits^n_{\R}\E^n\to \Q$  sending the cube $e_1\wedge\dots \wedge e_n\mapsto 1$, we see that $a_0\wedge\dots\wedge a_n$ will map to the volume of $\sigma$. So, the projection:
\[
\mathcal{P}(\E^n,T(n))\to \bigwedge\nolimits^n_{\R}\E^n\to\R
\]
is the volume map. There is a similar description when $q<n$. By definition, we see that  $x\in \mathcal{D}^q(\E^n)$ is an assignment, for each strict flag $\Phi=(U_0\supset\dots\supset U_{q})$ where $\mathrm{dim}(U_q)=q$, a volume element \[x_{\Phi}=v_1\wedge\dots\wedge v_q\in \bigwedge^q_\R U_q\] such that $\sum_{i=1}^q (-1)^i v_1\wedge\dots\wedge \hat{v_i}\wedge\dots\wedge v_q=0$. Given a non-degenerate simplex $\sigma=(a_0,\dots,a_n)$ we can associate an element $x(\sigma)\in \mathcal{D}^q$ as follows: we say $\sigma\ ||\  \Phi$ if there is some permutation $\pi\in \Sigma_n$ such that the affine space spanned by that $n-i+1$ simplex $(a_{\pi(i+1)},\dots,a_{\pi(n)})$ is parallel to $U_i$ for each $i=0,\dots,p$. We then define:
\[
x(\sigma)_{\Phi}=\begin{cases}
    0\quad&\text{ if $\sigma$ is not parallel to $\Phi$}\\
    \mathrm{sign}(\pi)(a_{\pi(n-q)}+\tau)\wedge\dots\wedge (a_{\pi(n)}+\tau)\quad&\text{ if $\sigma\ ||\  \Phi$}
\end{cases}
\]
where $\tau\in\E^n$ is the unique translation vector taking the affine space spanned by $(a_{\pi(i+1)},\dots,a_{\pi(n)})$ to the vector space $U_i$.
\end{remark}

\begin{corollary}\label{cor:Lusztigtrans}
    For $\E^n$, we can use the same strategy as in \ref{lustigex} to obtain following exact sequence of $O(n)$-modules:
    \[
    0\to \mathcal{D}^q(\E^n)\to \bigoplus_{U^{n-1}}\mathcal{D}^q(U^{n-1})\to \dots\to \bigoplus_{U^{q+1}}\mathcal{D}^q(U^{q+1})\to \bigwedge\nolimits_{\R}^q(U^{q})^t\to \bigwedge\nolimits_{\R}^q(\E^n)^t\to 0
    \]
    for $q<n$. Here, however, $U^{i}$ are now vector sub-spaces, and not affine. Here, the $t$ in the subscript stands for ``twist" as in \ref{twist}.
\end{corollary}
\begin{proof}
    Recall from \ref{lustigex} that we have a filtration of chain complexes on $C_*(V)$ for any affine subspace $V\subset \E^n$:
    \[
    \O=\Big(C_*(V)\Big)_{-1}\subset\Big(C_*(V)\Big)_{0}\subset \dots\subset \Big(C_*(V)\Big)_{n}=C_*(V)
    \]
    where $\Big(C_*(V)\Big)_{p}$ is defined as in \ref{simplicial chains}. This filtration induces a filtration on the co-invariants $(C_*(V))_{T(n)}$. Setting $V=\E^n$ we see that
    \[
    Gr_p((C_*(\E^n))_{T(n)})=\begin{cases}
        \bigoplus_{\mathrm{dim}(U_p)=p}C_*(U_p)_{T(U_p)}\quad & *=p\\
        0\quad &\text{otherwise}
    \end{cases}
    \]
    where $U_p$ now ranges over \emph{vector subspaces} of $\E^n$ as any two affine $k$-spaces are know equivalent by translations! Here, $T(U_k)$ denotes the translation group of $U_k$. As usual, the filtration gives rise to a spectral sequence with $E^1$ page: 
    \[
    E^{1}_{p,q}=H_{p+q}\Big(Gr_p((C_*(\E^n))_{T(n)})\Big)
    \] 
    Using our knowledge of $Gr_k(C_*(\E^n)_{T})$ and the proof of \ref{thm:sseqtr}, we have
    \[
    E^{1}_{p,q}=\begin{cases}
        \bigoplus_{\mathrm{dim}(U_p)=p}\bigoplus_{k=1}^p\mathcal{D}^k(U_p)\quad&q=0\\
        0\quad&\text{otherwise }
    \end{cases}
    \]
    further, the differential $E^1_{p,0}\to E^1_{p-1,0}$ is induced by the boundary map \[\bigoplus_{q=0}^p\partial_q:\bigoplus_{q=0}^pC_q(U_p)_{T(U_p)}\to \bigoplus_{q=0}^pC_{q-1}(U_p)_{T(U_p)}\]
    Since these are the only possible non-zero terms on the $E^1$ page, it follows that
    \[
    E^\infty_{p,q}=\begin{cases}
        \bigwedge\nolimits_\R^q(\E^n)\quad& p=0\\
        0\quad&\text{otherwise}
    \end{cases}
    \]
    
    Restricting ourselves to $\mathcal{D}^q(\E^n)\subset E^{0,n}$, we see that we have the sequence of maps:
    \[
    0\to \mathcal{D}^q(\E^n)\xrightarrow{\partial_q} \bigoplus_{U^{n-1}}\mathcal{D}^q(U^{n-1})\to\dots\to \bigoplus_{U^{q+1}}\mathcal{D}^q(U^{q+1})\to \bigoplus_{U^{q}}\mathcal{D}^q(U^q)\to \bigoplus_{U^{q-1}}\mathcal{D}^q(U^{q-1})=0
    \]
the last term is $0$ as there are no $q$-forms on a $q-1$ dimensional space. Recalling that $\partial$ is exact on $C_q(U_p)$ for $p>1$ like in \ref{lustigex}, it follows that it will be exact on $C_q(U_p)_{T(U_p)}$ for $p>0$ as well. At $p=0$, homology computes $\bigwedge^q_{\R}(\E^n)$ as desired. 
    
\end{proof}

\section{The Case of $\E^3$: First Theorems}
The goal of this section is to prove \ref{DehnSydlerJessen} by reducing the problem to \ref{thm:exceptional}, and then proving this. The idea here is that dimension $n=3$ is small enough to explicitly compute group homology. First we can do better than \ref{cor:decomp}.
\begin{corollary}\label{centerkills}
If $n+q\not \equiv 0 \ (\mathrm{mod}\ 2)$ then $H_0(O(n),(\mathcal{D}^q(\E^n)))=0$, and we have an isomorphism $H_0(O(n),\mathcal{D}^n(\E^n))\cong\R$ given by $n$-dimensional volume. Further:
\[
\mathcal{P}(\E^n)\cong \bigoplus_{n+q \equiv 0 \ (\mathrm{mod}\ 2)}H_0(O(n),(\mathcal{D}^q(\E^n)))
\]

\noindent
In particular, for $n=2$:
\[
\mathcal{P}(\E^2)\cong H_0(O(2),\mathcal{D}^2(\E^2))\cong\R
\]
so $2$-dimensional scissors congruence is just given by area. In dimension $3$:
\[
\mathcal{P}(\E^3)\cong H_0(O(3),\mathcal{D}^1(\E^3))\oplus \R
\]
where the projection onto $H_0(O(3),\mathcal{D}^3(\E^3))\cong \R$ denotes the volume. 
\end{corollary}
\begin{proof}
In definition \ref{twist}, we recall that the action of $\phi\in O(3)$ on $\mathrm{Pt}(\E^n)$, and therefore $\mathcal{D}^q$ must record whether or not $\phi$ is orientation preserving or reversing. We see that $-\mathrm{id}$ acts by a factor of $(-1)^q$ on $\bigwedge^q(U_q)$ for any $q$-dimensional subspace $U_q\subset\E^n$. Further, $\mathrm{-id}$ is orientation preserving iff $n$ is even, and so the sign of $\mathrm{id}$ acting on $\mathcal{D}^q$ is given by $(-1)^n$. This means that the action $\mathcal{D}^q$ is $(-1)^{q+n}$. If $n+q\not \equiv 0 \ (\mathrm{mod}\ 2)$, then every element $\mathcal{D}^q$ is $2$-torsion, meaning it must vanish as it is an $\R$-vector space. The isomorphism
\[H_0(O(n),(\mathcal{D}^n(\E^n))\cong \R\]
arises from the fact that we already have a volume isomorphism $\mathcal{D}^n(\E^n)\cong \bigwedge\nolimits_\R\E^n\cong \R$ as discussed in the previous section, and $O(n)$ (acting with sign) acts trivially on $\bigwedge\nolimits_\R\E^n$. 
\end{proof}
\begin{corollary}
    The inclusion $\mathcal{P}(\E^2)\xrightarrow{\times [0,1]}\mathcal{P}(\E^3)$ induced by sending a polygon to its unit-height prism $P\mapsto P\times [0,1]$ gives us the splitting:
    \[
    \mathcal{P}(\E^3)\cong \mathcal{P}(\E^2)\oplus H_0(O(3),\mathcal{D}^1) 
    \]
\end{corollary}
\begin{proof}
    Given any $V\in\R$, one can see that the cube $[0,V]\times [0,1]\times [0,1]$ has volume $V$, and so we have a surjective group homomorphism:
    \[
    \R=\mathcal{P}(\E^2)\xrightarrow{\times [0,1]} \mathcal{P}(\E^3)\rightarrow H_0(O(3),\mathcal{D}^3(\E^3))\cong \R
    \]
    and so it must actually be an isomorphism onto the summand.
\end{proof}
One can also make the following simplification via the \emph{Shapiro lemma}, a standard lemma in group homology (see 6.2 of \cite{Brown} for a proof).
\begin{lemma}[Shapiro]\label{shapirolemma}
    Let $H\subset G$ be a subgroup. Then, if $M$ is an $H$-module:
    \[
    H_*(G,\Z[G]\otimes_{\Z[H]} M)\cong H_*(H,M)
    \]
\end{lemma}
Using this lemma, along with a few other tricks we can compute group homology in simple cases:
\begin{theorem}\label{O(3)E1}
\[H_i(O(3),\bigoplus_{\E^1\subset\E^2}(\E^1)^t)=\begin{cases}
    \R\otimes_\Z\bigwedge^i_\Q(\R/\Z)\quad&\text{ if $i=1$}\\
        0\quad&\text{for $i=0,2$}
    \end{cases}\]
    where the direct sum is over all vector subspaces of $\R^3$.
\end{theorem}
\begin{proof}
    First, $O(3)$ acts transitively on all subspaces $\E^1$ of $\E^3$. The stabilizer of any such subspace is $O(1)\times O(2)$, where $O(1)=\Z/2\Z$ flips the copy of $\E^1$, and elements of $O(2)$ rotate around the copy of $\E^1$. In other words, we find that $O(3)/(O(2)\times O(1))$ parameterizes all lines of $\E^3$, and therefore we can understand the action of $O(3)$ on the module via left multiplication on $O(3)$.
    The action of of the stabilizer on any copy of $\E^1$ gives the module $(\E^1)^t\otimes \Z^t$ over $O(1)\times O(2)$, where $O(1)$ acts on $\E^1$ in the usual way, and $\Z^t$ denotes the action of $O(2)$ on $\Z$ with a sign depending on whether the action is orientation preserving or reversing. As $O(3)$ modules, there is an isomorphism: 
    \[
    \Z\cong \Z[O(3)]\bigotimes_{\Z[O(1)\times O(2)]}\Z
    \]
    and so 
    \[\bigoplus_{\E^1}(\E^1)^t\cong\Z[O(3)]\bigotimes_{\Z[O(1)\times O(2)]}\Big((\E^1)^t\otimes(\E^1)^t\Big)\]
where $O(1)$ and $O(2)$ act on separate copies of $(\E^1)^t$. However, $O(2)$ denotes rotations around $(\E^1)^t$ which means that $O(2)$ just acts by its determinant on $(\E^1)^t$. Therefore, as $O(2)$ modules $(\E^1)^t\cong \Z^t\otimes \R$, where $\Z^t$ denotes the action of $O(2)$ on $\Z$ by its determinant. 

By \ref{shapirolemma} one concludes:
\[
H_i(O(3),\bigoplus_{\E^1}(\E^1)^t)\cong \R\otimes_{\Z} H_i(O(1)\times O(2),(\E^1)^t\otimes\Z^t)
\]
By K\"{u}nneth, we have:
\[
H_*(O(1)\times O(2),(\E^1)^t\otimes(O(2))^t)\cong H_*(O(1),(\E^1)^t)\otimes_\Z H_*(O(2),\Z^t)
\]
As the action of $O(1)$ on $(\E^1)^t$ is trivial, $H_*(O(1),(\E^1)^t)=H_*(O(1))\otimes \R$. This will vanish in all nonzero degrees as $O(1)=\Z/2\Z$ is torsion. So:
\[
H_*(O(1)\times O(2),(\E^1)^t\otimes\Z^t)\cong  \R\otimes_\Z H_*(O(2),\Z^t)
\]
Next, observe that $O(2)=\{\mathrm{id},-\mathrm{id}\}\ltimes SO(2)$, where $SO(2)$ is the circle group. As a subgroup of $O(2)$, $SO(2)$ acts trivially on $\Z^t$ as all its matrices are determinant $1$, so we can check that we have an isomorphism of chain complexes $O(2)$-modules:
\[
C_*(O(2),\Z^t)=C_*(SO(2),\Z)^{\mathrm{det}}
\]
where the `$\mathrm{det}$' denotes the action of $g\in O(2)$ on the right-hand side sends a matrix $M\mapsto \mathrm{det}(g)\cdot M$. As $\R$ has no torsion:
\[
\R\otimes_\Z C_*(SO(2),\Z)\cong \R\otimes_\Z C_*(SO(2)/\mathrm{torsion})
\]
By the same argument as in \ref{thm:abtorfree} we have that \[H_n(SO(2)/\mathrm{torsion})\cong \bigwedge\nolimits_\Z^n (SO(2)/\mathrm{torsion})=\bigwedge\nolimits_\Q^n (SO(2)/\mathrm{torsion})\] 
Here $SO(2)/\mathrm{torsion}\cong \R/\Q$ is the $\Q$-vector space corresponding to the circle group modulo rational angles. In addition, $\R/\Q\cong \R/\Z$, where a particular isomorphism choice is chosen to be $\theta\mapsto \theta/\pi$. Taking homology, we have
\[
\R\otimes_\Z H_*(O(2),\Z^t)\cong \bigwedge\nolimits_{\R}^*(\R/\Z)^{\mathrm{det}}
\]
where $\Z/2\Z$ acts on $\R/\Z$ by negation $\theta\mapsto -\theta$. Combined with the determinant sign, $\Z/2\Z$ acts trivially on $\bigwedge\nolimits_{\R}^n(\R/\Z)$ for $n$ odd, and by $(-1)$ when $n$ is even. This means that the even groups are both $2$-torsion and $\R$-vector spaces, making them vanish, giving us the result:
\[
\R\otimes H_n(O(2),(\E^1)^t)\cong\begin{cases}
    \R\otimes \bigwedge\nolimits_{\R}^n(\R/\Z)\quad&\text{$n$ odd}\\
    0\quad&\text{$n$ odd}
\end{cases}
\]

\end{proof}
\begin{theorem}\label{thm:almostDSJ}
There is an exact sequence of groups:
\[
    0\to H_2(O(3),(\E^3)^t)\to H_0(O(3),\mathcal{D}^1(\E^3))\to \R\otimes\R/\Z\to H_1(O(3),(\E^3)^t)\to 0
    \]
\end{theorem}
\begin{proof}
    Setting $n=3$ and $q=1$ in corollary \ref{cor:Lusztigtrans}, one has the following exact sequence of $O(n)$-modules
    \[
    0\to \mathcal{D}^1(\E^3)\to \bigoplus_{\E^2\subset\E^3}\mathcal{D}^1(\E^2)\to \bigoplus_{\E^1\subset\E^2}\E^1\to (\E^3)^t\to 0
    \]
    Any exact sequence of modules gives rise to a spectral sequence on group homology with a vanishing $E_\infty$ page. The $E_1$ page is:
    \[\begin{tikzpicture}
  \matrix (m) [matrix of math nodes,
    nodes in empty cells,nodes={minimum width=5ex,
    minimum height=5ex,outer sep=-5pt},
    column sep=1ex,row sep=1ex]{
          2      &     H_2(O(3),(\E^3)^t)\   & \ H_2(O(3),\bigoplus_{\E_1\subset\E_2}(\E^1)^t) \ &  \ H_2(O(3),\bigoplus_{\E_2\subset\E_3}\mathcal{D}^1(\E^2)) \ & \ H_2(O(3),\mathcal{D}^1(\E^3)) \\
          1     &   H_1(O(3),(\E^3)^t) \  & \ H_1(O(3),\bigoplus_{\E_1\subset\E_2}(\E^1)^t) \ &  \ H_1(O(3),\bigoplus_{\E_2\subset\E_3}\mathcal{D}^1(\E^2)) \ & \ H_1(O(3),\mathcal{D}^1(\E^3))\\
          0     & H_0(O(3),(\E^3)^t) \  & \ H_0(O(3),\bigoplus_{\E_1\subset\E_2}(\E^1)^t) \ &  \ H_0(O(3),\bigoplus_{\E_2\subset\E_3}\mathcal{D}^1(\E^2))\  & \ H_0(O(3),\mathcal{D}^1(\E^3))\\
    \quad\strut &   0  &  1  &  2  & 3& \strut \\};
  \draw[-stealth] (m-3-3.west) -- (m-3-2.east);
  \draw[-stealth] (m-3-4.west) -- (m-3-3.east);
  \draw[-stealth] (m-3-5.west) -- (m-3-4.east);
    \draw[-stealth] (m-2-3.west) -- (m-2-2.east);
  \draw[-stealth] (m-2-4.west) -- (m-2-3.east);
  \draw[-stealth] (m-2-5.west) -- (m-2-4.east);
    \draw[-stealth] (m-1-3.west) -- (m-1-2.east);
  \draw[-stealth] (m-1-4.west) -- (m-1-3.east);
  \draw[-stealth] (m-1-5.west) -- (m-1-4.east);
\draw[thick] (m-1-1.east) -- (m-4-1.east) ;
\draw[thick] (m-4-1.north) -- (m-4-6.north) ;
\end{tikzpicture}\]

    By arguments similar to \ref{centerkills}, the second column completely vanishes as $1$ and $3$ are both odd. Combining this with our computations for the first column in \ref{O(3)E1}, the $E_1$ page is:

    \[\begin{tikzpicture}
  \matrix (m) [matrix of math nodes,
    nodes in empty cells,nodes={minimum width=5ex,
    minimum height=5ex,outer sep=-5pt},
    column sep=1ex,row sep=1ex]{
          2      &     H_2(O(3),(\E^3)^t)\   & \ 0 \ &  \ 0 \ & \ H_2(O(3),\mathcal{D}^1(\E^3)) \\
          1     &   H_1(O(3),(\E^3)^t) \  & \ \R\otimes\R/\Z \ &  \ 0 \ & \ H_1(O(3),\mathcal{D}^1(\E^3))\\
          0     & H_0(O(3),(\E^3)^t) \  & \ 0\ &  \ 0\  & \ H_0(O(3),\mathcal{D}^1(\E^3))\\
    \quad\strut &   0  &  1  &  2  & 3& \strut \\};
  \draw[-stealth] (m-3-3.west) -- (m-3-2.east);
  \draw[-stealth] (m-3-4.west) -- (m-3-3.east);
  \draw[-stealth] (m-3-5.west) -- (m-3-4.east);
  \draw[-stealth] (m-2-3.west) -- (m-2-2.east) ;
  \draw[-stealth] (m-2-4.west) -- (m-2-3.east);
  \draw[-stealth] (m-2-5.west) -- (m-2-4.east);
    \draw[-stealth] (m-1-3.west) -- (m-1-2.east);
  \draw[-stealth] (m-1-4.west) -- (m-1-3.east);
  \draw[-stealth] (m-1-5.west) -- (m-1-4.east);
\draw[thick] (m-1-1.east) -- (m-4-1.east) ;
\draw[thick] (m-4-1.north) -- (m-4-6.north) ;
\end{tikzpicture}\]
So, the $E_2$ page becomes
    \[\begin{tikzpicture}
  \matrix (m) [matrix of math nodes,
    nodes in empty cells,nodes={minimum width=5ex,
    minimum height=5ex,outer sep=-5pt},
    column sep=1ex,row sep=1ex]{
          2      &     H_2(O(3),(\E^3)^t)\   &\  0\  &\   0\  & \ H_2(O(3),\mathcal{D}^1(\E^3)) \\
          1     &   H_1(O(3),(\E^3)^t) \  & \ \mathrm{ker}(d^{1,1}_1)\ &  0 & \ H_1(O(3),\mathcal{D}^1(\E^3))\\
          0     & 0 \  & \ 0 \ &  0 & \ H_0(O(3),\mathcal{D}^1(\E^3))\\
    \quad\strut &   0  &  1  &  2  & 3& \strut \\};
\draw[-stealth] (m-3-5.north west) -- (m-2-3.south east);
\draw[thick] (m-1-1.east) -- (m-4-1.east) ;
\draw[thick] (m-4-1.north) -- (m-4-6.north) ;
\end{tikzpicture}\]
    Since the $E^\infty$ page must vanish, it follows that the image of $d^{3,0}_{2}$ must be $\mathrm{ker}(d^{1,1}_1)$ giving us the exact sequence:
    \begin{equation}\label{sequence:almost DSJ}
    H_0(O(3),\mathcal{D}^1(\E^3))\to \R\otimes\R/\Z\to H_1(O(3),(\E^3)^t)\to 0
    \end{equation}
    By doing a similar analysis for the $E^3$ page, we can complete \ref{sequence:almost DSJ} to the the exact sequence
    \[
    0\to H_2(O(3),(\E^3)^t)\to H_0(O(3),\mathcal{D}^1(\E^3))\to \R\otimes\R/\Z\to H_1(O(3),(\E^3)^t)\to 0
    \]
\end{proof}
\begin{corollary}
Recalling that $\mathcal{P}(\E^3)\cong \mathcal{P}(\E^2)\oplus H_0(O(3),\mathcal{D}^1(\E^3))$, we have a diagram of maps:
\[
\begin{tikzcd}
    0\arrow[r]&H_2(O(3),(\E^3)^t)\arrow[r]&H_0(O(3),\mathcal{D}^1(\E^3))\arrow[r]& \R\otimes\R/\Z\arrow[r]&H_1(O(3),(\E^3)^t)\arrow[r]&0\\
    0\arrow[r]&\mathcal{P}(\E^2)\arrow[r,"\times {[}0{,}1{]}"]&\mathcal{P}(\E^3)\arrow[u,two heads]\arrow[r,"D"]&\R\otimes\R/\Z\arrow[u,equal]&&
\end{tikzcd}
\]
where $D$ is the Dehn invariant, defined on $\mathcal{P}(\E^3)$ by sending a polytope $P$ with edges $\ell_i$ and dihedral angles $\theta_i$ to $D(P):=\sum_i\ell_i\otimes \theta_i/\pi$.
\end{corollary}
\begin{proof}
The proof of this theorem will revolve around understanding the geometric information encoded in the $E_1$ differential. Reproducing the exact sequence $n=3$ and $q=1$ in corollary \ref{cor:Lusztigtrans}, we have:
\[
    0\to \mathcal{D}^1(\E^3)\xrightarrow{\partial_1^{\E^3}} \bigoplus_{\E^2\subset\E^3}\mathcal{D}^1(\E^2)\xrightarrow{\partial_1^{\E^2}} \bigoplus_{\E^1\subset\E^2}(\E^1)^t\xrightarrow{\partial_1^{\E^1}}  (\E^3)^t\to 0
\]
 we can consider the short exact sequences
\[
    0\to \mathcal{D}^1(\E^3)\to \bigoplus_{\E^2\subset\E^3}\mathcal{D}^1(\E^2)\to \mathrm{ker}(\partial_1^{\E^1})\to 0
\]
this gives us a long exact sequence in homology:
\[\dots\to H_1(O(3),\bigoplus_{\E^2\subset\E^3}\mathcal{D}^1(\E^2))\to H_1(O(3),\mathrm{ker}(\partial_1^{\E^1}))\xrightarrow{\delta_1} H_0(O(3),\mathcal{D}^1(\E^3))\to H_0(O(3),\bigoplus_{\E^2\subset\E^3}\mathcal{D}^1(\E^2))\]
By the same arguments as in \ref{centerkills}, we see that $H_*(O(2),\bigoplus_{\E^2\subset\E^3}\mathcal{D}^1(\E^2))=0$ because $2$ and $1$ do not have the same parity. This gives us an isomorphism $H_1(O(3),\mathrm{ker}(\partial_1^{\E^1}))\cong H_0(O(3),\mathcal{D}^1(\E^3))$.

Consider the simplex $\sigma=|(0,v_1,v_2,v_3)|$, and denote $[x(\sigma)]\in H_0(O(3),\mathcal{D}^1(\E^3)$ as in the discussion after \ref{cor:decomp}. We see that
\[
x(\sigma)_{\R\l v_1,v_2\r\supset\R\l v_1\r}=v_1
\]
and 
\[
x(\sigma)_{\R\l v_1,v_3\r\supset\R\l v_1\r}=-v_1
\]
Consider the element:
\[
(\mathrm{id},\theta_{v_1})\otimes v_1\in C_1(O(3),\bigoplus_{\E_1\subset \E^2}(\E^1)^t)
\]
where $\theta_{v_1})\in SO(2)$ is the rotation matrix around $v_1$ by the dihedral angle of $v_1$. Then one can check on homology:
\[
\delta_1([(\mathrm{id},\theta_{v_1})\otimes v_1)])_{U_0\supset U_1}=\begin{cases}
    v_1\quad& U_0=\R\l v_1,v_2\r, U_1=\R\l v_1\r\\
    -v_1\quad& U_0=\R\l v_1,v_3\r, U_1=\R\l v_1\r\\
    0\quad&\text{otherwise}
\end{cases}
\]
Now, the map:
\[
H_0(O(3),\mathcal{D}^1(\E^3))\xrightarrow{\cong} H_1(O(3),\mathrm{ker}(\partial^{\E^1}_1))\to H_1(O(3),\bigoplus_{\E^1\subset\E^2}(\E^1)^t)
\]
is exactly the differential $d_2^{3,0}$, and the isomorphism $H_1(O(3),\bigoplus_{\E^1\subset\E^2}(\E^1)^t)\cong \R\otimes_\Z\R/\Z$ sends $[(\mathrm{id},\theta)\otimes v_1)]$ to $||v_1||\otimes \theta$, with $||v_1||$ is the length of $v_1$. Doing this for the other vertices of $\sigma$, we get that
$[x(\sigma)]$ is identified with $\sum_{i}[(1,\theta_{v_i})\otimes ||v_i||]\in H_1(O(3),\bigoplus_{\E^1\subset\E^2}(\E^1)^t)$, where $\theta_{v_i}$ is the rotation around $v_i$ (translated to the origin) by dihedral angle of $v_i$, and $i$ runs over all vertices. This gives the full Dehn invariant.
\end{proof}
\section{The Case of $\E^3$: Exceptional Isomorphisms}
So now, all that is left is the following two theorems:
\begin{theorem}
    The Dehn invariant is injective, i.e. $H_2(O(3),(\E^3)^t)\cong 0$.
\end{theorem}

\begin{theorem}
    The cokernel of Dehn invariant is  $H_1(O(3),(\E^3)^t)\cong \Omega^{1}_{\R}$, and the induced map
    \[
    \R\otimes\R/\Z\to \Omega^{1}_{\R}
    \]
    from \ref{thm:almostDSJ} sends $\ell\otimes \theta\mapsto \ell\frac{d\cos(\theta)}{\cos{(\theta)}}$
\end{theorem}
First, we prove a small lemma:
\begin{lemma}\label{subduetorsion}
    Let $G=\Gamma\ltimes N$, where $\Gamma$ is a torsion group. Then, for a torsion-free $G$-module $A$ such that $\Gamma$ acts trivially on $A$, we have:
    \[
    H_*(G,A)\cong H_*(N,A)_{\Gamma}
    \]
    where $H_*(N,A)_{\Gamma}$ are the coinvariants.
\end{lemma}
\begin{proof}
    This follows from the Hochshild-Serre spectral sequence: The $E^2$ is given by:
    \[
    E^2_{p,q}=H_p(\Gamma,H_q(N,A))\Rightarrow H_{p+q}(G,A)
    \]
    However, $H_q(N,A)$ will be torsion-free as $A$ is torsion-free, and so $ E^2_{p,q}=0$ unless $p=0$. Therefore, the $E^\infty$ page reads:
    \[
    H_{n}(G,A)=H_0(\Gamma, H_n(N,A))=H_n(N,A)_{\Gamma}
    \]
\end{proof}
\begin{definition}\label{def:Hochschild}
    Given an algebra $A$ over $\Q$, define $\epsilon_i^n:A^{\otimes (n+1)}\to A^{\otimes n}$ such that \[\epsilon_i^n(a_0\otimes\dots\otimes a_n)=a_0\otimes\dots\otimes a_ia_{i+1}\otimes\dots\otimes a_n\] Define
    \[
    \Omega_n(A):=\bigcap_{0\leq i\leq n-1}\mathrm{ker}(\epsilon_i^n)
    \]
    there is a boundary map $b_n:\Omega_n(A)\to \Omega_{n-1}(A)$ given by:
    \[
    b_n(a_0\otimes\dots\otimes a_n)=\sum_{0\leq i\leq n-1}(-1)^ia_0\otimes\dots\otimes a_ia_{i+1}\otimes\dots\otimes a_n+(-1)a_na_0\otimes a_1\otimes\dots\otimes a_{n-1}
    \]
    One can check that $b_{n-1}\circ b_n=0$, making $(\Omega_n(A),b)$ a chain complex. Write $I_n(A)=\mathrm{ker}(b_n)$ and $B_n(A)=\mathrm{im}(b_{n+1})$. The \emph{Hochshild Homology} of $A$ is defined to be\begin{footnote}{The original definition of Hochshild Homology different. For our purposes, we can just take this as the definition. See \cite{Weibel} Section 8.8.1 and 9.4.2 for details.}
    \end{footnote}:
    \[
    HH(A)=I_n(A)/B_n(A)
    \]
\end{definition}

\begin{remark}\label{remark:differentials}
    For $a\in A$, the \emph{differential} of $a$ is defined as $da:=1\otimes a-a\otimes 1$. One can show that all elements of $\Omega_n(A)$ are $\Q$-linear combinations of the form $a_0da_1\dots da_n$, where $a_i\in A$. Observe that $d(a_i+k)=da_i$ for $k\in\Z$. Further, there is an isomorphism
    \[
    \Omega_n(A)\cong A\otimes (A/\Q)^{\otimes n}
    \]
    mapping $a_0da_1\dots da_n\mapsto a_0\otimes a_1\otimes\dots\otimes a_n$.
\end{remark}
\begin{definition}
  Given a commutative ring $R$, let $\Omega_R^1$ denote the ring of K\"{a}hler differntials of $R$ (over $\Z$). The generators are $dr$, for every $r\in R$ and $dn=0$ for $n\in\Z$. Further, for $r,s\in R$, there are two relations:
  \[
  d(r+s)=dr+ds;\quad\quad d(rs)=r(ds)+s(dr)
  \]
  Further, define $\Omega^n_{R}:=\bigwedge_{\Q}\Omega^1_{R}$. 
\end{definition}
The following is proven in 9.4.7 in \cite{Weibel}, and for $*=1$ (which is what is most important for us) one can check the statement by hand.
\begin{theorem}[Hochschild-Kostant-Rosenberg]\label{HKR}
    If $\F$ is a field of characteristic $0$, then 
    \[
    HH_*(K)\cong \Omega^*_K
    \]
    
\end{theorem}

\begin{theorem}\label{O(n)Pin(n)}
    Let $\H$ denote the quaternion algebra, $\mathrm{Pin}(4)$ denote the double cover of $O(4)$, and $(\H\otimes\H)^{-}$ denote the negative eigenspace of the conjugation map $\tau:a_0\otimes a_1\mapsto a_0^*\otimes a_1^*$. As $\tau$ commutes with $b$ defined in \ref{def:Hochschild}, we can define $I_1(\H)^{-}\subset (\H\otimes\H)^{-}$ to the $(-1)$-eigenspace of $\tau$ acting on $I_1(\H)$. Then $I_1(\H)^{-}$ can be given the structure of a $\mathrm{Pin}(4)$ module, and we have isomorphisms for $i=1,2$:
    \[H_i(O(3),(\E^3)^t)\cong H_{i-1}(\mathrm{Pin}(4),I_1(\H)^{-})\] 
\end{theorem}
\begin{proof}
    First, we have $H_i(O(3),(\E^3)^t)\cong H_i(SO(3),\E^3)$ by \ref{subduetorsion} and the action of $\pm id$ on $(\E^3)^t$ is trivial due to the sign. Next, let $\mathrm{Spin}(n)$ be the double cover of $SO(n)$, and $Sp(1)\subset\H$ denote the unit quaternions. We state a few exceptional isomorphisms:
    \begin{enumerate}
        \item $\mathrm{Spin}(3)\cong \mathrm{Sp}(1)$
        \item $\mathrm{Spin}(4)\cong \mathrm{Sp}(1)\times \mathrm{Sp}(1)$
        \item $\mathrm{Pin}(4)\cong \Z/2\Z\ltimes (\mathrm{Spin}(3)\times \mathrm{Spin}(3))$, where $\Z/2\Z$ interchanges the copies of $\mathrm{Spin}(3)$.
    \end{enumerate}
     For $\sigma^+=(1, q_1,q_2)\in \mathrm{Pin}(4)$ where $q_i\in \mathrm{Sp}(1)$, define its action on $(\H\otimes\H)^{-}$ as \[\sigma^+(a_1\otimes a_2)=q_1a_1q_2^*\otimes q_2a_2q_1^*\] for $a_1\otimes a_2\in\H\otimes\H$. If instead we have $\sigma^-=(-1, q_1,q_2)$, then define $\sigma^-(a_1,a_2):=\sigma^+(a_2^*,a_1^*)$. Similarly, $\mathrm{Pin}(4)$ has an action on the direct sum $\H\oplus \H$, defined, for $\sigma^+=(1, q_1,q_2)\in \mathrm{Pin}(4)$, by conjugation in each variable:
    \[
    \sigma(a_1,a_2)=(q_1a_1q_1^*,q_2a_2q_2^*)
    \]
    for $(a_1,a_2)\in\H\oplus \H$, and $\sigma^-(a_1,a_2):=\sigma^+(-a_2^*,-a_1^*)$. By definition of $I_1(\H)$ in \ref{def:Hochschild}, we have a long exact sequence of $\mathrm{Pin}(4)$-modules:
    \[
    0\to I_1(\H)^{-}\to (\H\otimes\H)^-\xrightarrow{(\epsilon_0,-\epsilon_1)} \H^-\oplus\H^-\to 0
    \]
    where the $\epsilon_i$ were defined in \ref{def:Hochschild}. This gives us a long exact sequence in group homology:
    \begin{equation}\label{lespin}
    \dots\to H_i(\mathrm{Pin}(4),I_1(\H)^{-})\to H_i(\mathrm{Pin}(4),(\H\otimes\H)^-)\to H_i(\mathrm{Pin}(4),\H^-\oplus\H^-)\to \dots
    \end{equation}
    As $G=\mathrm{Pin}(4)=\Z/2\Z\ltimes \mathrm{Spin}(4)$, and $\H^-\oplus\H^-$ is torsion-free, \ref{subduetorsion} gives us:
    \begin{equation}\label{equationpin}
    H_n(\mathrm{Pin}(4),\H^-\oplus\H^-)= H_{n}(\mathrm{Spin}(4),\H^-\oplus\H^-)_{\Z/2\Z}
    \end{equation}
    We try to now understand $H_{n}(\mathrm{Spin}(4),\H^-)=H_{n}(\mathrm{Spin}(4),\H^-)\oplus H_{n}(\mathrm{Spin}(4),\H^-)$. By K\"{u}nneth:
    \begin{equation}\label{hard}
    H_{n}(\mathrm{Spin}(4),\H^-)=\bigoplus_{p+q=n}H_p(\mathrm{Spin}(3),\H^-)\otimes_{\H^-} \Big(H_q(\mathrm{Spin}(3))\otimes_\Z \H^-\Big)
    \end{equation}
    We see that $H_0(\mathrm{Spin}(3))=\Z$, and $H_1(\mathrm{Spin}(3))$ is its abelianization which is $0$ as it is a perfect group. Lastly, by an argument similar to \ref{centerkills}, we have $H_2(\mathrm{Spin}(3))$ is $2$-torsion, implying $H_2(\mathrm{Spin}(3))\otimes_\Z\H^-=0$. Putting all of this together, we get for $n=0,1,2$:
    \[
    H_{n}(\mathrm{Spin}(4),\H^-)=H_{n}(\mathrm{Spin}(3),\H^-)
    \]
    Further $\H^-=\R\l \hat{i},\hat{j},\hat{k}\r$ and $\mathrm{Spin}(4)=\Z/2\Z\ltimes SO(3)$. We see that that $\Z/2\Z$ acts trivially on $\H^-$, and $SO(3)$ acts by rotations, and so by \ref{subduetorsion}, for all $n$:
    \[
    H_{n}(\mathrm{Spin}(3),\H^-)=H_{n}(SO(3),\R^3)
    \]
   Now, the action of $\Z/2\Z$ in \ref{equationpin} simply flips the copies of $\mathrm{Spin(3)}$ and also $\H^-$ by exceptional isomorphism (3), giving us for $i=0,1,2$
   \[
   H_{i}(\mathrm{Pin}(4),\H^-\oplus\H^-)=H_i(SO(3),\R^3)
   \]
In the appendix we show that $H_i(\mathrm{Pin}(4),(\H\otimes\H)^-)=0$ for $i=0,1,2$, and so the long exact sequence in \ref{lespin} gives us the isomorphisms for $i=1,2$:
    \[H_i(SO(3),\R^3)\cong H_{i-1}(\mathrm{Pin}(4),I_1(\H)^{-})\] 
\end{proof}
Before we can finish the proof, we need another few definitions and lemmas:
\begin{lemma}\label{hochschildspin}
$(\Omega_*(\H),b)$ can be given the structure of a chain complex of $\mathrm{Spin}(4)$ modules, with conjugation action extended from $\H=\Omega_0$. Further, $\mathrm{Spin}(4)$ acts trivially on $HH_n(\H)$ for all $n\in\N$
\end{lemma}
\begin{proof}
    Define, for $a_0\in\H$ and for $\sigma=(q_1,q_2)\in\mathrm{Spin}(4)$,  $\sigma(a_0)=q_2a_0q_2^*$, and for $\H^n$ for $n>0$, let
    \[
    \sigma(a_0\otimes\dots\otimes a_n)=q_1a_0q_2^*\otimes q_2a_1q_1^*\otimes q_1a_2q_1^*\dots\otimes q_1a_2q_1^*
    \]
    One can compute that if $\epsilon_i^n(a_1\otimes\dots\otimes a_n)=0$, then $\epsilon_i^n(\sigma(a_1\otimes\dots\otimes a_n))=0$ as well. One can also check that $\H\otimes\H/\Q =\mathrm{ker}(\epsilon_0^1)=\Omega_1(\H)$, and so the isomorphism in \ref{remark:differentials} can be promoted to an isomorphism of $\mathrm{Spin}(4)$-modules:
    \begin{equation}\label{isomorphismspin4mod}
    \Omega_n(\H)\cong \Omega_1(\H)\otimes (\H/\Q)^{\otimes (n-1)}
    \end{equation}
    where the action on the right side is given, for $\sigma=(q_1,q_2)$ by \[\sigma(\omega\otimes da_2\dots da_n)=\sigma(\omega)\otimes d(q_1a_2q_1^*)\otimes \dots \otimes d(q_na_2q_n^*)\]
    for $\omega\in\Omega_1(\H)$
    Further, the conjugation action $\tau$ on $\H=\Omega_0$ extends to one on $\H\otimes\H$, and one can check that on $\Omega^1(\H)$, $\tau(a_0(da_1))=-a_0^*d(a_1^*)$. If $\omega=a_0da_1$ then write $\omega^*=a_0^*da_1^*$. Using \ref{isomorphismspin4mod}, we define conjugation on $\Omega_n$ by:
    \[
    \tau(\omega\otimes d_2\otimes\dots\otimes da_n)=(-1)^{(n-2)(n-3)/2} \omega^*\otimes  da^*_n\otimes\dots\otimes da^*_2
    \]
    $\tau$ commutes with the differential, making $\pm 1$ eigenspaces of $\Omega_*$, denoted $\Omega_*^\pm$, chain complexes themselves. To see that $\mathrm{Spin}(4)$ acts trivially on $HH_n(\H)$, for any $(q,1)\in Spin(4)$, pick a field $\F$ containing it and $\R$. Then, for $a\in\F$ conjugation by $q$ is trivial because $\F$ is commutative, and it will act trivially on $HH_n(\F)$ which implies it will also act trivially on $HH_n(\H)$. For $(1,q)$, the proof is the same, and so all elements of $\mathrm{Spin}(4)$ are generated by these, concluding the proof.
\end{proof}
\begin{lemma}\label{lem:spin(3)action}
    For any of the two summands $\mathrm{Spin}(3)\subset \mathrm{Spin}(4)$ one has $H_0(\mathrm{Spin}(3),\H)=\R$, and
    for $n>0$:
\[
H_0(\mathrm{Spin}(3),\Omega_n(\H))=0
\]
\end{lemma}
\begin{proof}
    For $n=0$, we see that for $q\in\mathrm{Spin}(3)$ and $a_0\in\H$, $(q,1)(a_0)=qa_0$, and so $H_0(\mathrm{Spin}(3),\H)=\H/\mathrm{Sp}(1)=\R$. For $n>0$, the isomorphism of $\mathrm{Spin}(4)$-modules in \ref{isomorphismspin4mod} tells us that we only need to prove that $H_0(\mathrm{Spin}(3),\Omega_1(\H))=0$, as everything lies in the same orbit as $0\otimes da_2\otimes\dots\otimes da_n=0$. Looking at the exact sequence of modules $0\to\Omega_1(\H)\to \H\otimes\H\xrightarrow{\epsilon_0}\H\to 0$ where $\mathrm{Spin}(3)$ acts trivially on the last module, we recall that $H_1(\mathrm{Spin}(3),\Z)=0$ as it is a perfect group, so we have the exact sequence:
    \[
    0\to H_0(\mathrm{Spin}(3),\Omega_1(\H))\to H_0(\mathrm{Spin}(3),\H\otimes\H)\xrightarrow{(\epsilon_0)_*} \H\to  0
    \]
    where $(\epsilon_0)_*(a_0\otimes a_1)=a_0a_1$. Consider the map $\eta:\H\to H_0(\mathrm{Spin}(3),\H\otimes\H)$ given by $\eta(a_0)=a_0\otimes 1$. The proof concludes if we show $\eta$ is surjective, i.e. every element of $ H_0(\mathrm{Spin}(3),\H\otimes\H)$ is of the form $a\otimes 1$. To see this, first note that it holds for $q\otimes a$ for $q\in\mathrm{Sp}(1)$, and so by addition, $(q+q^*)\otimes a$. Now, $[-2,2]\subset \R\subset\H$ can be written as $q+q^*\in\R$ for some $q$. So, by addition, the same holds for real numbers $r\otimes a$, and since all $x\in\H$ are of the form $rq$ where $r\in\R,q\in \mathrm{Sp}(1)$, we see that it holds always.
\end{proof}
\begin{lemma}
    We have an isomorphism \[HH_n(\H)\cong HH_n(\R)=\Omega^1_\R\] 
\end{lemma}
\begin{proof}
    Let $\H(\Q)\subset\H$ denote the $\Q$-quaternion algebra. Then $\H(\Q)\otimes\R=\H$, and so by K\"{u}nneth:
    \[
    HH_*(\H)=HH_*(\H(\Q))\otimes HH_*(\R)
    \]
    We show that $HH_0(\H(\Q))=\Q$ and $HH_n(\H(\Q))=0$ for $n>0$, and the conclusion follows. There is a map of $\R$-algebras $\H\otimes\H\to \mathrm{End}_\R(\H)\cong M_4(\R)$ sending $a_0\otimes a_1$ to the function $v\mapsto a_0va_1^{-1}$. The kernel is $0$ as $\H\otimes\H$ is a simple algebra, and checking the dimensions, they are the same. So,  $\H(\Q)\otimes\H(\Q)\cong M_4(\Q)$. Now, by 9.5.3 of \cite{Weibel}, $M_4(\Q)$ and $\Q$ are \emph{Morita invariant}, i.e. their categories of left (and right) modules are equivalent, and so by 9.5.6 of \cite{Weibel} (it is more general than what we need), we have $HH_*(M_4(\Q))=HH_*(\Q)=\Omega^*_{\Q}$. So, again by K\"{u}nneth:
    \begin{align*}
      HH_*(\H(\Q))\otimes HH_*(\H(\Q))&\cong HH_*(\H(\Q)\otimes \H(\Q))\\
      &\cong HH_*(M_4(\Q))\\
      &HH_*(\Q)\cong \Omega^*_{\Q}  
    \end{align*}  
    Just by checking the definition we see that $\Omega^n_{\Q/\Z}=0$ for $n=1$ (and therefore $n>1$) and $\Q$ for $n=0$. This forces $HH_0(\H(\Q))=\Q$ and $HH_n(\H(\Q))=0$ for $n>0$, concluding the proof.
\end{proof}
\begin{theorem}
    The quotient map $I_1(\H)\to HH_1(\H)$ induces an isomorphism:
    \[
    H_0(\mathrm{Pin}(4),I_1(\H)^-)\cong HH_1(\H)\cong \Omega^1_{\R}
    \]
\end{theorem}
\begin{proof}
    For the first part of the theorem, we see that extending $(\Omega_*,b)$ to a chain complex of $\mathrm{Spin}(4)$ modules compatible with conjugation $\tau$ gives us the exact sequence of $(-1)$ $\tau$-eigenspaces:
    \[
    0\to B_1(\H)^-\to I_1(\H)^-\to HH_1(\H)^-\to 0
    \]
    one can compute that the involution of $\tau$ on $HH_1(\H)$ is given by $-1$, just by computing it for $\R\subset\H$, as this inclusion induces an isomorphism on Hochshild Homology. So $HH_1(\H)^-=HH_1(\H)$, and since $\mathrm{Spin}(4)$ acts trivially on $HH_1(\H)$ we get a long exact sequence:
    \[
    H_1(\mathrm{Spin}(4))\otimes HH_1(\H)\to H_0(\mathrm{Spin}(4),B_1(\H)^-) \to H_0(\mathrm{Spin}(4),I_1(\H)^-) \to HH_1(\H) \to 0 
    \]
    One sees that $H_1(\mathrm{Spin}(4))=0$ as its abelianization is zero, and $H_0(\mathrm{Spin}(4),B_1(\H)^-)=0$ due to the long exact sequence in homology of $I_2^-\to (\Omega_2)^-\to B_1\to 0$ and by K\"{u}nneth: 
    \[
    H_0(\mathrm{Spin}(4),\Omega_2^-)=H_0(\mathrm{Spin}(3),\Omega_2^-)\otimes H_0(\mathrm{Spin}(3),\Omega_2^-)=0
    \]
    due to \ref{lem:spin(3)action}. Finally, since $\mathrm{Pin}(4)\cong\Z/2\Z\ltimes \mathrm{Spin}(4)$, and the copy of $\Z/2\Z$ acts trivially on all of $(\Omega_*,b)$, by \ref{subduetorsion} we get an induced isomorphism $H_0(\mathrm{Pin}(4),I_1(\H)^-)\cong HH_1(\H)=\Omega^1_\R$. 
    \end{proof}
    \begin{theorem}
    The induced map
    \[
    \R\otimes_{\Z}\R/\Z\to H_1(SO(3),\E^3)\xrightarrow{\cong} H_0(\mathrm{Pin}(4),I_1(\H)^-)\cong\Omega^1_\R
    \]
    sends $\ell\otimes\theta\mapsto \ell\frac{d\cos{(\theta)}}{\sin{(\theta)}}$
    \end{theorem}
    \begin{proof}
    Recall from \ref{thm:almostDSJ} that $J$ is the composition 
    \[
    \R\otimes H_1(O(2),\Z^t)\to H_1(O(3),(\E^3)^t)=H_0(\mathrm{Pin}(4),HH_1(\H))
    \]
    where $\Z^t$ is the $O(2)$-module corresponding to rotations of about a line in $\E^3$ (up to sign). For the purposes of this proof let $\E^1$ denote the line about with $O(2)$ rotates - then $H_1(O(2),\Z^t)=H_1(O(2),(\E^1)^t)$. Since both $O(2)=\Z/2\Z\ltimes O(2)$ and $\mathrm{Pin}(2)=\Z/2\Z\ltimes U(1)$, the double cover $\rho:\mathrm{Pin}(2)\to O(2)$ sends $(\pm 1,e^{i\theta})\mapsto (\pm 1, e^{2i\theta})$. There action of $\mathrm{Pin}(2)$ on the plane (via rotation and conjugation) perpendicular to $i\R\subset\C$. This defines a $\mathrm{Pin}(2)$ module $i\R$, and as $\rho$ is an isomorphism, we see that the action on $i\R$ is the same as the action of $O(2)$ on $(\E^1)^t$, giving us an isomorphism $H_*(O(2),(\E^1)^t)\cong H_*(\mathrm{Pin}(2),i\R)$. Now, the usual inclusion $\C\hookrightarrow\H$ commutes with conjugations, and so we have a commutative diagram:
    \[
    \begin{tikzcd}
    0\arrow[r]& I_1(\C)^-\arrow[d]\arrow[r] &(\C\otimes\C)^-\arrow[r,"\epsilon_0"]\arrow[d]&  i\R\arrow[d,"\Delta"]\arrow[r] &0\\
    0\arrow[r]& I_1(\H)^-\arrow[r] &(\H\otimes\H)^-\arrow[r,"\epsilon"]&  \H^-\otimes\H^-\arrow[r]& 0
    \end{tikzcd}
    \]
where $\Delta(ir)=(ir,-ir)$ for $ir\in i\R$. We can embed $Pin(2)\hookrightarrow Pin(4)$ by sending $(\pm 1, e^{i\theta})\mapsto (\pm 1, e^{i\theta},e^{-i\theta})$, where we view $e^{i\theta}\in U(1)\subset \mathrm{Sp}(1)$ as a unit quaternion. So, we get a commutative diagram:
\begin{equation}\label{pin2pin4}
\begin{tikzcd}
    H_1(\mathrm{Pin}(2),i\R)\arrow[r]\arrow[d,"\Delta_*"]& H_0(\mathrm{Pin}(2),I_1(\C)^-)\arrow[d]\arrow[r] &H_0(\mathrm{Pin}(2),HH_1(\C)^-)\arrow[d,"\cong"]\\
    H_1(\mathrm{Pin}(4),\H^-\otimes\H^-)\arrow[r]& H_0(\mathrm{Pin}(4), I_1(\H)^-)\arrow[r,"\cong"] & H_0(\mathrm{Pin}(4), HH_1(\H))
    \end{tikzcd}
\end{equation}
where the right-most horizontal morphisms are induced by the module maps $\C\otimes\C\to \mathrm{ker}(\epsilon_0)/B_1(\C)$. Since $\C^-=i\R$, we have $HH_1(\C)^-=HH_1(\C^-)=HH_1(\R)=\Omega_1^\R$. We showed the following string of isomorphisms \[H_0(\mathrm{Pin}(4), I_1(\H)^-)\cong H_0(\mathrm{Pin}(4), HH_1(\H))\cong \Omega_1^\R\] 
and in \ref{O(n)Pin(n)}, we showed that $H_1(\mathrm{Pin}(4),\H^-\otimes\H^-)\cong H_1(SO(3),\E^3)$ so the map $2\cdot J$ is now the map $H_1(\mathrm{Pin}(2),i\R)\xrightarrow{\Delta_*} H_1(\mathrm{Pin}(4),\H^-\otimes\H^-)$, where the factor of $2$ comes from the double cover $\rho$. We therefore just need to study the top row in \ref{pin2pin4} in order to understand $J$, and most importantly, the connecting homomorphism in the top row of \ref{pin2pin4}. In order to do this, first identify $(\ell\otimes\theta/2\pi)\in \R\otimes\R/\Z$ with $(e^{i\theta/2}\otimes i\ell)$ in $C_1(\mathrm{Pin}(2),i\R)$. Next, studying $\epsilon_0$ in \ref{pin2pin4}, we see that $\epsilon_0$ sends $e^{i\theta/2}\otimes (i\ell \otimes 1+1\otimes i\ell)\in C_1(\mathrm{Pin}(2),(\C\otimes\C)^-)$ to $e^{i\theta/2}\otimes 2i\ell$. So, the connecting homomorphism in \ref{pin2pin4} sends $e^{i\theta/2}\otimes 2i\ell$ to the image of $e^{i\theta/2}\otimes (i\ell \otimes 1+1\otimes i\ell)$ in the boundary map $C_1(\mathrm{Pin}(2),(\C\otimes\C)^-)\to C_0(\mathrm{Pin}(2),(\C\otimes\C)^-)$, i.e. to 
\[
(i\ell \otimes 1+1\otimes i\ell)-e^{i\theta/2}=i\ell \otimes 1+1\otimes i\ell-ie^{i\theta}\otimes e^{-i\theta}-e^{i\theta}\otimes ie^{-i\theta}
\]
which in $\Omega^1(\C)^-$ is $\frac{i\ell }{2}(e^{-i\theta}de^{i\theta}-e^{i\theta} de^{i\theta})$. Using Euler's formula $e^{i\theta}=i\sin(\theta)+\cos(\theta)$, one can compute
\begin{align*}  
e^{-i \theta}de^{i\theta}&=[\cos(\theta)-i\sin(\theta)]\ d[\cos(\theta)+i\sin(\theta)+]\\
&=\cos(\theta)d\cos(\theta)+\sin(\theta)d\sin(\theta)+i[\cos(\theta)d\sin(\theta)-\sin(\theta)d\cos(\theta)]\\
&=i[\cos(\theta)d\sin(\theta)-\sin(\theta)d\cos(\theta)]
\end{align*}
    as $\cos(\theta)d\cos(\theta)+\sin(\theta)d\sin(\theta)=\frac{1}{2}d[\cos^2(\theta)+\sin^2(\theta)]=0$. This identity allows us to replace $d\sin(\theta)$ with $\frac{\cos(\theta)}{\sin(\theta)}d\cos(\theta)$. Doing so in the last step, we get
    \begin{align*}
        e^{-i \theta}de^{i\theta}&=i\Big[\frac{-\cos^2(\theta)-\sin^2(\theta)}{\sin(\theta)}d\cos(\theta)\Big]\\
        &=\frac{-i}{\sin(\theta)}d\cos(\theta)
    \end{align*}
    This also implies $e^{-i \theta}de^{i\theta}=\frac{-i}{\sin(\theta)}d\cos(\theta)$, giving us that 
    \[
    \frac{i\ell }{2}(e^{-i\theta}de^{i\theta}-e^{i\theta} de^{i\theta})=\frac{\ell}{\sin(\theta)}d\cos(\theta) 
    \]
     as desired
\end{proof}
All that needs to be proved is:
\begin{theorem}
    \[H_1(\mathrm{Pin}(4),I_1(\H)^-)=0\]
\end{theorem}
\begin{proof}
Due to \ref{hochschildspin},both $I_2(\H)$ and $B_2(\H)$ are $\mathrm{Spin}(4)$ modules with a conjugation action. One can see that for $\R\subset \C$, $\tau$ acts by identity. The induced morphism $HH_*(\R)\to HH_*(\C)$ commutes with conjugation and is an isomorphism, so one has $HH_*(\C)^-=0$, implying that $I_2(\H)^-=B_2(\H)^-$. Then, we have the exact sequence $\mathrm{Spin}(4)$-modules:
\begin{equation}\label{vanishses1}
0\to I_3(\H)^-\to \Omega_3(\H)^-\xrightarrow{b_3} I_2(\H)^-\to 0
\end{equation}
\begin{equation}\label{vanishses2}
0\to I_2(\H)^-\to \Omega_2(\H)^-\xrightarrow{b_2} I_1(\H)^-\to 0
\end{equation}
we know that since $H_0(\mathrm{Spin}(4),\Omega_3(\H)^-)=0$ by \ref{lem:spin(3)action} and the K\"{u}nneth formula, we have that $H_0(\mathrm{Spin}(4),I_2(\H)^-)=0$ by the long exact sequence in homology arising from \ref{vanishses1}. So putting this into \ref{vanishses2}, we have \[H_1(\mathrm{Spin}(4),\Omega_2(\H)^-)\xrightarrow{(b_1)_*}H_1(\mathrm{Spin}(4),I_1(\H)^-)\to 0\]
so $(b_1)_*$ is surjective. The result follows if we prove that $(b_1)_*=0$. For simplicity, let $S_\pm$ denote the two copies of $\mathrm{Spin}(3)$. By the Hochschild-Serre spectral sequence for $\mathrm{Spin}(4)=S_+\times S_-$ yields the usual 5-term exact sequence:
\[
H_2(S_+,H_0(S_-,\Omega_2(\H)^-))\to H_0(S_+,H_1(S_-,\Omega_2(\H)^-))\to H_1(\mathrm{Spin}(4),\Omega_2(\H)^-)\to H_1(S_+,H_0(S_-,\Omega_2(\H)^-))\to 0
\]

From \ref{lem:spin(3)action}, we have $H_0(S_-,\Omega_2(\H)^-)=H_0(S_-,\Omega_2(\H)^-)=0$ giving the isomorphism $H_0(S_+,H_1(S_-,\Omega_2(\H)^-))\cong H_1(\mathrm{Spin}(4),\Omega_2(\H)^-)$. One therefore has a commutative diagram:
\[\begin{tikzcd}
    H_0(S_+,H_1(S_-,\Omega_2(\H)^-))\arrow[r,"(b_1)_*"]&H_1(\mathrm{Spin}(4),I_1(\H)^-)\\    H_1(S_-,\Omega_2(\H)^-)\arrow[u]\arrow[r]&H_1(S_-,I_1(\H))\arrow[u]
\end{tikzcd}\]
where the map $H_1(S_-,\Omega_2(\H)^-)\to H_0(S_+,H_1(S_-,\Omega_2(\H)))$ is given by quotienting out the $S_+$ action. Remembering the module structure on $\Omega_2(\H)=\Omega_1(\H)\otimes (\H/\Q)$ in \ref{hochschildspin}, we can choose the copy of $S_-$ that does not act on the second factor, making $H_1(S_-,\Omega_2(\H)^-)=H_1(S_-,I_1(\H)^-)\otimes(\H/\Q)^-$. What we therefore need to is that the image of $H_1(S_-,I_1(\H)^-)$ must be zero in $H_1(\mathrm{Spin}(4),I_1(\H)^0)$. Elements of $H_1(S_-,I_1(\H)^-)$ are of the form $\xi\cdot d(i\alpha),\xi\cdot d(i\alpha)$ or $\xi\cdot d(k\alpha)$ where $\alpha\in\R$, as $\H^-=\R\l \hat{i},\hat{j},\hat{k}\r$. Suppose, without loss of generality, it is $\xi\cdot d(i\alpha)$. Then, we show that taking $S_+$ orbits: $\xi\cdot d(i\alpha)=0\in H_0(H_1(S_-,I_1(\H)^-))$. If $q=x+iy\in U(1)\subset S_+$, we have the action of $q$ on $\xi$ is conjugation, i.e. $(x+iy)\xi (x-iy)$. Similarly one can subtract this expression to the action of $q^*$ on $\xi$ to get that, modulo the action of $S_+$, 
\[
2i(x\xi y -y\xi x)=0\in H_0(S_+,H_1(S_-,I_1(\H)^-))
\]
As $I_1(\H)$ is an $\R$-module, we multiply on the left by $\frac{1}{2y}$ and on the right by $\frac{1}{y}$, and setting $\alpha=\frac{x}{y}$, we get $i\alpha\xi-\xi i\alpha=0$, or $\xi\  d(i\alpha)=0$, where we can take $x,y$ to achieve any $\alpha\in\R$.
\end{proof}
\appendix
\section{The group homology of $O(n)$ and $E(n)$}
The purpose of this appendix is to prove that $H_i(\mathrm{Pin}(4),(\H\otimes\H)^-)=0$ for $i=0,1,2$. The proof, however, follows from more general results from in chapter 9 of \cite{Dupont} that are interesting in their own right. We first rephrase the vanishing theorem we must prove. There is an isomorphism of $\mathrm{Pin}(4)$-modules $(\H\otimes\H)^-\cong \bigwedge_{\R}^2(\H)$, with isomorphism sending $a_0\otimes a_1\mapsto a_0\wedge a_1^*$. Viewing $\H=\R^4$ as an $\R$ vector space, we see that the action of $\mathrm{Pin}(4)$ in \ref{O(n)Pin(n)} is compatible with the $O(4)$ action $\R^4$. Further, $\mathrm{Pin}(4)=\Z/2\Z\ltimes O(4)$, where $\Z/2\Z$ changes the action of $O(4)$ on $\H=\R^4$ by a sign. As a result, the sign change will be trivial on $\bigwedge_{\Z}^2(\H)$, and so we find by \ref{subduetorsion}:
\[
H_*(\mathrm{Pin}(4),(\H\otimes\H)^-)=H_*(O(4),\bigwedge\nolimits^2_{\R}(\R^4))
\]
So, we must now show that the right-hand side vanishes for $*=0,1,2$. 
\begin{theorem}
For all $n\in\N$, there is an isomorphism
    \[
    H_k(\mathrm{Isom}(\E^n))\cong H_k(O(n))\oplus\bigoplus_{q=1}^k H_{k-q}\Big(O(n),\bigwedge\nolimits_\Q^q(\R^n)\Big)
    \]
\end{theorem}
\begin{proof}
    The Hochshild-Serre spectral sequence for $\mathrm{Isom}(\E^n)=O(n)\ltimes T(n)$ gives us the $E_2$ page:
    \[
    E^2_{p,q}=H_{p}(O(n),H_q(T(n)))\Rightarrow H_{p+q}(\mathrm{Isom}(\E^n),\Z)
    \]
We have already computed in \ref{thm:abtorfree} that since $T(n)\cong \R^n$, we have $H_q(T(n))=\bigwedge\nolimits^q_\R(\R^n)$ for $q>0$ and $H_0(T(n))=\Z$. By the same argument as in \ref{thm:sseqtr}, the $E^2$ differentials vanish due to the the fact that they commute with the dilations $\mu_n:\alpha\mapsto n\alpha$. The result follows.
\end{proof}
Combined with the previous theorem, the vanishing in low degrees follows from:
\begin{theorem}
    The inclusion $O(n)\hookrightarrow \mathrm{Isom}(\E^n)$ induces isomorphisms: $H_k(O(n))\to H_k(\mathrm{Isom}(\E^n))$
    for $k<n$.
\end{theorem}
\begin{proof}
    The fact that $H_k(O(n))\to H_k(\mathrm{Isom}(\E^n))$ is injective follows from the split injection of chain complexes $C_*(O(n))_{O(n)}\to C_*(\mathrm{Isom}(\E^n))_{\mathrm{Isom}(\E^n)}$ (i.e. $H_k(O(n))\to H_k(\mathrm{Isom}(\E^n))$ is injective for all $k$). We recall that the complex in \ref{simplicial chains} denoted $C_*(\E^n)$ is acylic, and inherits a natural action of $\mathrm{Isom}(\E^n)$, so one has a spectral sequence with $E^1$ page:
    \[
    E^1_{p,q}=H_p(\mathrm{Isom}(\E^n),C_q(\E^n))\Rightarrow H_{p+q}(\mathrm{Isom}(\E^n),\Z)
    \]
    called the \emph{hyperhomology spectral sequence} (the more well-known dual notion is the hypercohomology spectral sequence). Setting $q=0$, $C_0(\E^n)=\E^n$, and $\mathrm{Isom}(\E^n)$ acts on $\E^n$ transitively with stabilizer $O(n)$, giving us an isomorphism of $\mathrm{Isom}(\E^n)$ modules:
    \[
    \E^n\cong \Z[\mathrm{Isom}(\E^n)]\otimes_{O(n)}\Z
    \]
    We see that by the Shapiro lemma \ref{shapirolemma}: \[E^1_{p,0}=H_p(\mathrm{Isom}(\E^n),\E^n)=H_p(O(n))\] 
    and for $0<q<n$ one can inductively show that $E^1_{p,q}=H_p(\mathrm{Isom}(\E^{n-q}))\otimes C_p(\E^q)_{\mathrm{Isom}(\E^q)}$ via the Shapiro lemma as well. The differential $d^1_{p,q}:E^1_{p,q}\to E^1_{p-1,q}$ is induced by the boundary map:
    \[
    C_p(\E^q)_{\mathrm{Isom}(\E^q)}\to C_{p-1}(\E^q)_{\mathrm{Isom}(\E^q)}
    \]
    We show that this boundary map is acyclic for $p+q\leq n$. Consider the complex $\mathrm{Gr}_*(C_*(\E^n))\subset C^*(\E^n)$ defined in \ref{prop:chains}. The simplicies of $\mathrm{Gr}_i(C_i(\E^n))$ are in \emph{general position}, i.e. for all $p\leq n$ there is no $(a_{i_0},\dots,a_{i_p})$ inside a $(p-1)$-dimensional affine subspace. We study the inclusion:
    \begin{equation}\label{laststand}
    H_k(\mathrm{Isom}(\E^n), \mathrm{Gr}_*(C_*(\E^n)))\to H_k(\mathrm{Isom}(\E^n),C_*(\E^n))
    \end{equation}
    Our work in \ref{prop:chains} implies that the induced map $\mathrm{Gr}_*(C_*(\E^n)))\to C_*(\E^n)$ is surjective, so we now show that this map is zero. A concrete formula for barycentric subdivision of a simplex $\sigma=(a_0,\dots,a_k)$ is defined inductively by:
    \[
    \mathrm{sd}(\sigma)=c_{\sigma}*\mathrm{sd}(\partial \sigma)
    \]
    where $c_{\sigma}$ is the barycenter and $*$ denotes the join operation. If $\sigma\in Gr_p(C_p(\E^n))$, then it also has a well-defined circumcenter inside of it, and is equal to the barycenter. The map for $p\leq n$ (otherwise there are no non-degenerate simplices):
    \[
    \mathrm{sd}_p:Gr_p(C_p(\E^n))\to C_p(\E^n)
    \]
    sends \[\mathrm{sd}(\sigma)=\sum_{\sigma=\sigma_0\supset \dots\supset\sigma_k}\pm (c_{\sigma_0},\dots,c_{\sigma_k})\]
    where the sign is positive if the induced orientation on $(c_{\sigma_0},\dots,c_{\sigma_k})$ is the same as $\sigma$, and negative otherwise. However, this sum is zero - if $\sigma_k=(a_1)$ and $\sigma_{k-1}=(a_1,a_2)$, for example, then they will cancel with the simplex indexed by $\sigma_k'=(a_2)$ and $\sigma_{k-1}=(a_1,a_2)$ because they have different sign, and the isometry that flips $(a_1,a_2)$ identifies $(a_1)$ with $(a_2)$. Therefore, $\mathrm{sd}_p=0$, and by the same proof as in \ref{prop:chains}, we have that for $k\leq n$, barycentric subdivision is chain homotopic to the inclusion map. This shows that \ref{laststand} is $0$ on homology for $k\leq n$.

    This shows that $d^1_{p,q}$ is acyclic, we find:
    \[
    E^2_{p,q}=\begin{cases}
        H_p(O(n))\quad& q=0,p<n\\
        0\quad& 0<q\leq n-p\\
    \end{cases}
    \]
    This proves that $H_{p}(\mathrm{Isom}(\E^n))\cong H_p(O(n))$ for $p<n$.
\end{proof}

\bibliographystyle{alpha}
\bibliography{refs}

\end{document}